\RequirePackage{ifpdf}
\ifpdf 
\documentclass[pdftex]{sigma}
\else
\documentclass{sigma}
\fi

\begin{document}


\renewcommand{\thefootnote}{$\star$}

\renewcommand{\PaperNumber}{065}

\FirstPageHeading

\ShortArticleName{$\mathrm{\frak sl}(2)$-Trivial Deformations of
${\rm Vect_{Pol}}(\mathbb{R})$-Modules of Symbols}

\ArticleName{$\boldsymbol{\mathrm{\frak sl}(2)}$-Trivial
Deformations\\ of $\boldsymbol{{\rm
Vect_{Pol}}(\mathbb{R})}$-Modules of Symbols\footnote{This paper
is a contribution to the Special Issue on Deformation
Quantization. The full collection is available at
\href{http://www.emis.de/journals/SIGMA/Deformation_Quantization.html}{http://www.emis.de/journals/SIGMA/Deformation\_{}Quantization.html}}}

\Author{Mabrouk BEN AMMAR  and  Maha BOUJELBENE}

\AuthorNameForHeading{M.~Ben Ammar and M.~Boujelbene}

\Address{D\'epartement de Math\'ematiques, Facult\'e des Sciences
de Sfax, BP 802, 3038 Sfax, Tunisie}
\Email{\href{mailto:mabrouk.benammar@fss.rnu.tn}{mabrouk.benammar@fss.rnu.tn},
\href{mailto:maha.boujelben@fss.rnu.tn}{maha.boujelben@fss.rnu.tn}}

\ArticleDates{Received January 14, 2008, in f\/inal form September
05, 2008; Published online September 18, 2008}

\Abstract{We consider the action of ${\rm Vect_{Pol}}(\mathbb{R})$
by Lie derivative on the spaces  of symbols of dif\/ferential
operators. We study the deformations of this action that become
trivial once restricted to $\rm\frak{sl}(2)$. Necessary and
suf\/f\/icient conditions for integrability of inf\/initesimal
deformations are given.}


\Classification{17B56; 17B66; 53D55}

\section{Introduction}\label{sec1}

Let ${\rm Vect_{Pol}}(\mathbb{R})$ be the Lie algebra of
polynomial vector f\/ields on $\mathbb{R}$. Consider the
1-parameter action of ${\rm Vect_{Pol}}(\mathbb{R})$  on the space
$\mathbb{R}[x]$ of polynomial functions on $\mathbb{R}$ def\/ined
by
\begin{gather*} L_{X\frac{d}{dx}}^\lambda(f)= Xf'+\lambda
X'f,\end{gather*} where $X, f\in\mathbb{R}[x]$ and
$X':=\frac{dX}{dx}$. Denote by $\cal{ F}_\lambda$ the
$\mathrm{Vect}_\mathrm{Pol}(\mathbb{R})$-mo\-du\-le structure on
$\mathbb{R}[x]$ def\/ined by this action for a f\/ixed $\lambda$.
Geometrically, $\cal{ F}_\lambda$ is the space of polynomial
weighted densities of weight $\lambda$ on $\mathbb{R}$
\begin{gather*}{\cal F}_\lambda=\big\{ fdx^{\lambda}\mid f\in
\mathbb{R}[x]\big\}.
\end{gather*}
The space ${\cal F}_\lambda$ coincides with the space of vector
f\/ields, functions and dif\/ferential 1-forms for $\lambda = -1,
0$ and $1$, respectively.

Denote by ${\cal D}_{\nu,\mu}:=\mathrm{Hom}_{\text{dif\/f}}({\cal
F}_\nu, {\cal F}_\mu)$ the ${\rm Vect_{Pol}}(\mathbb{R})$-module
of linear dif\/ferential operators with the ${\rm
Vect_{Pol}}(\mathbb{R})$-action given by the formula
\begin{gather}\label{Lieder2}L_X^{\nu,\mu}(A)=L_X^\mu\circ
A-A\circ L_X^\nu.\end{gather} Each module $\mathcal{D}_{\nu,\mu}$
has a natural f\/iltration by the order of dif\/ferential
operators; the graded mo\-du\-le ${\cal S}_{\nu,\mu}:=\mathrm{gr}
\mathcal{D}_{\nu,\mu}$ is called the {\it space of symbols}. The
quotient-module
$\mathcal{D}^k_{\nu,\mu}/\mathcal{D}^{k-1}_{\nu,\mu}$ is
isomorphic to the module of tensor densities
$\mathcal{F}_{\mu-\nu-k} $, the isomorphism is provided by the
principal symbol~$\sigma$ def\/ined by
\begin{gather*}
A=\sum_{i=0}^ka_i(x)\partial_x^i\mapsto\sigma(A)=a_k(x)(dx)^{\mu-\nu-k}
\end{gather*} (see, e.g., \cite{gmo}). As a ${\rm
Vect_{Pol}}(\mathbb{R})$-module, the space ${\cal S}_{\nu,\mu}$
depends only on the dif\/ference $\delta=\mu-\nu$, so that ${\cal
S}_{\nu,\mu}$ can be written as ${\cal S}_{\delta}$, and we have
\begin{gather*}
{\cal S}_{\delta} = \bigoplus_{k=0}^\infty \mathcal{F}_{\delta-k}
\end{gather*} as
${\rm Vect_{Pol}}(\mathbb{R})$-modules. The space of symbols of
order $\leq n$ is
\begin{gather*}
{\cal S}_\delta^n:=\bigoplus_{j=0}^n{\cal F}_{\delta-j}.
\end{gather*}

The space ${\cal D}_{\nu,\mu}$ cannot be isomorphic as a ${\rm
Vect_{Pol}}(\mathbb{R})$-module to the corresponding space of
symbols, but is a deformation of this space in the sense of
Richardson--Neijenhuis~\cite{nr2}; however, they are isomorphic as
$\frak{sl}(2)$-modules (see~\cite{hg}). In the last two decades,
deformations of various types of structures have assumed an ever
increasing role in mathematics and physics. For each such
deformation problem a goal is to determine if all related
deformation obstructions vanish and many beautiful techniques had
been developed to determine when this is so. Deformations of Lie
algebras with base and versal deformations were already considered
by Fialowski~\cite{f1}. It was further developed, with
introduction of a complete local algebra base (local means a
commutative algebra which has a unique maximal ideal) by
Fialowski~\cite{f2}. Also, in~\cite{f2}, the notion of miniversal
(or formal versal) deformation was introduced in general, and it
was proved that under some cohomology restrictions, a versal
deformation exists. Later Fialowski and Fuchs, using this
framework, gave a construction for the versal
deformation~\cite{ff2}.

We use the framework of Fialowski~\cite{f2} (see also~\cite{abbo}
and~\cite{aalo2}) and consider (multi-parameter) deformations over
complete local algebras. We construct the miniversal deformation
of this action and def\/ine the complete local algebra related to
this deformation.

According to Nijenhuis--Richardson~\cite{nr2}, deformation theory
of modules is closely related to the computation of cohomology.
More precisely, given a Lie algebra $\frak{g}$ and a
$\frak{g}$-module $V$, the inf\/initesimal deformations of the
$\frak{g}$-module structure on $V$, i.e., deformations that are
linear in the parameter of deformation, are related to
$\mathrm{H}^1\left(\mathfrak{g};\mathrm{End}(V)\right)$. The
obstructions to extension of any inf\/initesimal deformation to a
formal one are related to
$\mathrm{H}^2\left(\mathfrak{g};\mathrm{End}(V)\right)$. More
generally, if $\frak h$ is a subalgebra of $\frak g$, then the
$\frak h$-relative cohomology space
$\mathrm{H}^1\left(\mathfrak{g},\frak h;\mathrm{End}(V)\right)$
measures the inf\/initesimal deformations that become trivial once
the action is restricted to $\frak h$ ($\frak h$-{\it trivial
deformations}), while the obstructions to extension of any $\frak
h$-trivial inf\/initesimal deformation to a~formal one are related
to $\mathrm{H}^2\left(\mathfrak{g},\frak h;\mathrm{End}(V)\right)$
(see, e.g.,~\cite{b}).

Denote $\mathcal{D}:=\mathcal{D}(n,\delta)$ the ${\rm
Vect_{Pol}}(\mathbb{R})$-module of dif\/ferential operators on
${\cal S}_\delta^n$. The inf\/initesimal deformations of the ${\rm
Vect_{Pol}}(\mathbb{R})$-module ${\cal S}_\delta^n$ are
classif\/ied by the f\/irst dif\/ferential cohomology space,
\[
\mathrm{H}^1_{\rm dif\/f}({\rm Vect_{Pol}}(\mathbb{R});{\cal
D})=\oplus_{\lambda,k}\mathrm{H}^1_{\rm dif\/f}({\rm
Vect_{Pol}}(\mathbb{R});{\cal D}_{\lambda,\lambda+k}),
\]
 while the
obstructions for integrability of inf\/initesimal deformations
belong to the second di\-f\/fe\-ren\-tial cohomology space,
\[
\mathrm{H}^2_{\rm dif\/f}({\rm Vect_{Pol}}(\mathbb{R});{\cal
D})=\oplus_{\lambda,k}\mathrm{H}^2_{\rm dif\/f}({\rm
Vect_{Pol}}(\mathbb{R});{\cal D}_{\lambda,\lambda+k}),
\] where,
hereafter,  $\delta-\lambda$ and $k$ are integers satisfying
$\delta-n\leq\lambda,\lambda+k\leq\delta$.

Here we study the ${\rm \frak{sl}}(2)$-trivial deformations, thus
we consider the dif\/ferential $\frak{sl}(2)$-relative cohomology
spaces. The f\/irst space
\[
{\rm H}^1_{\rm dif\/f}({\rm Vect_{Pol}}(\mathbb{R}),{\rm
\frak{sl}}(2);{\cal D})=\oplus_{\lambda,k}{\mathrm H}^1_{\rm
dif\/f}({\rm Vect_{Pol}}(\mathbb{R}),{\rm \frak{sl}}(2) ;{\cal
D}_{\lambda,\lambda+k})
\] was calculated by Bouarroudj and
Ovsienko~\cite{bo}; and  the
second space
\[
{\mathrm H}^2_{\rm dif\/f}({\rm
Vect_{Pol}}(\mathbb{R}),\frak{sl}(2);{\cal
D})=\oplus_{\lambda,k}{\mathrm H}^2_{\rm dif\/f}({\rm
Vect_{Pol}}(\mathbb{R}),\frak{sl}(2);{\cal D}_{\lambda,\lambda+k})
\]was calculated by Bouarroudj \cite{b}. We give explicit expressions of some 2-cocycles that span the
cohomology group ${\mathrm H}^2({\rm Vect_P}(\mathbb{R}),{\rm
\frak{sl}}(2) ;{\cal D}_{\lambda,\lambda+k})$.

This paper is organized as follows. In Section~\ref{sec2} we study
some properties of the $\frak{sl}(2)$-invariant dif\/ferential
operators. These properties are related to the
$\frak{sl}(2)$-relative cohomology. In Section~\ref{sec3} we study
the f\/irst and the second $\frak{sl}(2)$-relative cohomology
spaces which are closely related to the deformation theory.
Especially we explain some $\frak{sl}(2)$-relative 2-cocycles
which naturally appear as obstructions to integrate any
$\frak{sl}(2)$-trivial inf\/initesimal deformation to a formal
one. In Section~\ref{sec4} we give an outline of the general
deformation theory: def\/initions, equivalence, integrability
conditions and miniversal deformations. In Section~\ref{sec5} we
give the f\/irst main result of this paper: Theorem \ref{Main1}.
That is, we explain all second-order integrability conditions for
any inf\/initesimal $\frak{sl}(2)$-trivial deformation of the
${\rm Vect_{Pol}}(\mathbb{R})$-module  ${\cal S}_\delta^n$. In
Section~\ref{sec6} we complete the list of integrability
conditions by computing those of third and fourth-order. We prove that these
conditions are necessary and suf\/f\/icient to integrate any
inf\/initesimal $\frak{sl}(2)$-trivial deformation to a formal
one. Moreover, we prove that any $\frak{sl}(2)$-trivial
deformation is, in fact, equivalent to a polynomial one of degree
$\leq2$: Theorem \ref{Main2}. Finally, in Section~\ref{sec7}, we
complete our study by giving a few examples of deformations.

\section[Invariant differential operators]{Invariant dif\/ferential operators}\label{sec2}

In this paper we study the ${\rm \frak{sl}}(2)$-trivial
deformations of the space of symbols of dif\/ferential operators
which is a ${\rm Vect_{Pol}}(\mathbb{R}) $-module, so we begin by
studying some properties of ${\rm \frak{sl}}(2)$-invariant
bilinear dif\/ferential operators.

 Let us consider the space of bilinear dif\/ferential
operators
$c:\mathcal{F}_\lambda\times\mathcal{F}_\mu\rightarrow\mathcal{F}_{\tau}$.
The Lie algebra, ${\rm Vect_{Pol}}(\mathbb{R})$, acts on this
space by the Lie derivative:
\begin{gather*}
L_X(c)(fdx^\lambda,\,gdx^\mu)=L_X^\tau\big(c(fdx^\lambda,gdx^\mu)\big)-c\big(L_X^\lambda
(fdx^\lambda),\,gdx^\mu\big)-c\big( fdx^\lambda,\,L_X^\mu
(gdx^\mu)\big).
\end{gather*}
A bilinear dif\/ferential operator
$c:\mathcal{F}_\tau\times\mathcal{F}_\lambda\rightarrow\mathcal{F}_{\mu}$
is called ${\rm \frak{sl}}(2)$-invariant if, for all $X\in{\rm
\frak{sl}}(2)$, we have $L_X(c)=0$, or equivalently
\begin{gather}\label{bilinvar1}
L_X^\mu\big(c(fdx^\tau,\,gdx^\lambda)\big)=c\big(L_X^\tau
(fdx^\tau),\,gdx^\lambda\big)+c\big( fdx^\tau,\,L_X^\lambda
(gdx^\lambda)\big).
\end{gather}
That is, the set of such ${\rm \frak{sl}}(2)$-invariant bilinear
dif\/ferential operators is the subspace on which the subalgebra
${\rm \frak{sl}}(2)$ acts trivially.

Now, let us consider a linear map $c:{\rm
Vect_{Pol}}(\mathbb{R})\rightarrow\mathcal{D}_{\lambda,\mu}$, then
we can see $c$ as a bilinear dif\/ferential operator $c:{\rm
Vect_{Pol}}(\mathbb{R})\times\mathcal{F}_\lambda\rightarrow\mathcal{F}_{\mu}$
(or
$c:\mathcal{F}_{-1}\times\mathcal{F}_\lambda\rightarrow\mathcal{F}_{\mu}$
since ${\rm Vect_{Pol}}(\mathbb{R})$ is isomorphic to
$\mathcal{F}_{-1}$) def\/ined by
$c(X,fdx^\lambda)=c(X)(fdx^\lambda).$ So, the $\frak{
sl}(2)$-invariance property (\ref{bilinvar1}) of $c$ reads, for
all $X\in{\rm \frak{sl}}(2)$ and $Y\in{\rm
Vect_{Pol}}(\mathbb{R})$,
\begin{gather*}
L_X^{\mu}\circ c(Y)=c([X,Y])+c(Y)\circ L_X^{\lambda},
\end{gather*}
or equivalently
\[
L_X^{\lambda,\mu}(c(Y))=c([X,Y]).
\]

The $\mathfrak{sl}(2)$-invariant bilinear dif\/ferential operators
were calculated by Gordan. We recall here  the results and we need
to add some precision concerning the space of the
$\mathfrak{sl}(2)$-invariant dif\/ferential operators from ${\rm
Vect_{Pol}}(\mathbb{R})$ to $\mathcal{D}_{\lambda,\mu}$ vanishing
on $\mathfrak{sl}(2)$.

\begin{proposition}[\cite{pg}]\label{trans2}  There exist
$\frak{sl}(2)$-invariant bilinear differential operators, called
transvectants,
\begin{gather*}
J_k^{\tau,\lambda}: \ \
\mathcal{F}_\tau\times\mathcal{F}_\lambda\rightarrow\mathcal{F}_{\tau+\lambda+k},\qquad
(\varphi dx^\tau,\phi dx^\lambda)\mapsto
J_k^{\tau,\lambda}(\varphi,\phi)dx^{\tau+\lambda+k}
\end{gather*} given by
\begin{gather*}
J_k^{\tau,\lambda}(\varphi,\phi)=\sum_{i+j=k}c_{i,j}\varphi^{(i)}\phi^{(j)},
\end{gather*} where $k\in\mathbb{N}$ and the coefficients $c_{i,j}$ are characterized
as follows:
\begin{itemize}\itemsep=0pt
  \item [\rm i)] If neither $\tau$ nor $\lambda$ belong to the set
$\{0,-{1\over2},-1,\dots,-{k-1\over2}\}$ then
\[
c_{i,j}=(-1)^{j}
\begin{pmatrix}2\tau+k\\j\end{pmatrix}
\begin{pmatrix}2\lambda+k\\i\end{pmatrix},
\] where  $\big(^x_i\big)$ is the standard binomial coefficient
 $\big(^x_i\big)=\frac {x(x-1)\cdots (x-i+1)}{i!}$.
  \item [\rm ii)] If  $\tau$ or $\lambda\in\{0,-{1\over2},-1,\dots,-{k-1\over2}\}$,
the coefficients $c_{i,j}$  satisfy the recurrence relation
\begin{gather}\label{cij}(i+1)(i + 2\tau)c_{i+1,j} + (j+1)(j+2\lambda)c_{i,j+1} = 0.\end{gather}
Moreover, the space of solutions of the system~\eqref{cij} is
two-dimensional if $2\lambda=-s$ and $2\tau=-t$ with $t > k-s-2$,
and one-dimensional otherwise.

\item[\rm iii)] For $k\geq3$, the space of
$\mathfrak{sl}(2)$-invariant linear differential operators from
${\rm Vect_{Pol}}(\mathbb{R})$ to ${\cal D}_{\lambda,\lambda+k-1}$
vanishing on $\mathfrak{sl}(2)$ is one-dimensional.
\end{itemize}
\end{proposition}

\begin{proof} We need to prove only part iii), for the other statements see, for instance, \cite{pg} or~\cite{b}. First, we recall that ${\rm Vect_{Pol}}(\mathbb{R})$ is
isomorphic to ${\cal F}_{-1}$ as ${\rm
Vect_{Pol}}(\mathbb{R})$-module. So, according to the formulae
(\ref{cij}), if $k\geq3$, the space of
$\mathfrak{sl}(2)$-invariant bilinear dif\/ferential operator from
${\rm Vect_{Pol}}(\mathbb{R})\times\cal{F}_\lambda$ to
$\mathcal{F}_{\lambda-1+k}$ is 2-dimensional if and only if
$2\lambda\in\{1-k,2-k,3-k\}$. Let us consider the transvectant
$J_k^{-1,\lambda}$ def\/ined by, for $\varphi{d\over dx}\in{\rm
Vect_{Pol}}(\mathbb{R})$ and $\phi dx^\lambda\in\cal{F}_\lambda$,
\begin{gather}\label{J}
J_k^{-1,\lambda}(\varphi,\phi)=\sum_{i+j=k,i\geq 3}
c_{i,j}\varphi^{(i)}\phi^{(j)},
\end{gather}where  the coef\/f\/icients $c_{i,j}$
satisfy (\ref{cij}). If $2\lambda\in\{1-k,2-k,3-k\}$ the space of
$\mathfrak{sl}(2)$-invariant linear dif\/ferential operator from
${\rm Vect_{Pol}}(\mathbb{R})$ to ${\cal D}_{\lambda,\lambda+k-1}$
is spanned by $J_k^{-1,\lambda}$ and $I_k^{-1,\lambda}$ where
\begin{gather*}I_k^{-1,\lambda}(\varphi,\phi)=
\left\{\begin{array}{lll}\varphi\phi^{(k)}\quad&\text{if}\quad 2\lambda=1-k,\\
\varphi\phi^{(k)}+{k\over2}\varphi'\phi^{(k-1)}\quad&\text{if}\quad
2\lambda=2-k,\\
\varphi\phi^{(k)}+k\varphi'\phi^{(k-1)}+{k(k-1)\over2}\varphi''\phi^{(k-2)}\quad&\text{if}\quad
2\lambda=3-k.
\end{array}\right.\end{gather*}
If $2\lambda\notin\{1-k,2-k,3-k\}$ the corresponding space is
one-dimensional and it is spanned by~$J_k^{-1,\lambda}$. We see
obviously that only the operators $J_k^{-1,\lambda}$ vanish on
$\frak{sl}(2)$. Part~iii) of Proposition~\ref{trans2} is proved.
 \end{proof}

\section{Cohomology spaces}\label{sec3}

Let $\frak{g}$ be a Lie algebra acting on a space $V$ and let
$\frak h$ be a subalgebra of $\frak g$. The space of $\frak
h$-relative $n$-cochains of $\frak{g}$ with values in $V$ is the
$\frak{g}$-module
\[
C^n(\frak{g},\frak{h}; V ) := \mathrm{Hom}_{\frak
h}(\wedge^n(\frak{g}/\frak{h});V).
\] The
{\it coboundary operator} $ \partial^n: C^n(\frak{g},\frak{h}; V
)\rightarrow C^{n+1}(\frak{g},\frak{h}; V )$ is a $\frak{g}$-map
satisfying $\partial^n\circ\partial^{n-1}=0$. The kernel of
$\partial^n$, denoted $Z^n(\mathfrak{g},\frak{h};V)$, is the space
of $\frak h$-relative $n$-{\it cocycles}, among them, the elements
in the range of $\partial^{n-1}$ are called $\frak h$-relative
$n$-{\it coboundaries}. We denote $B^n(\mathfrak{g},\frak{h};V)$
the space of $n$-coboundaries.

By def\/inition, the $n^{th}$ $\frak h$-relative  cohomolgy space
is the quotient space
\[
H^n
(\mathfrak{g},\frak{h};V)=Z^n(\mathfrak{g},\frak{h};V)/B^n(\mathfrak{g},\frak{h};V).
\]
We will only need the formula of $\partial^n$ (which will be
simply denoted $\partial$) in degrees 0, 1 and 2: for $v \in
C^0(\frak{g},\frak{h}; V ) = V^{\frak h}$, $\partial v(X) := Xv$,
for $ b\in C^1(\frak{g}, \frak{h};V )$,
\[
\partial
b(X, Y ) := Xb(Y)-Y b(X) -b([X, Y ])
\]
and for $\Omega\in C^2(\frak{g}; \frak{h},V )$,
\begin{gather}\label{2c}
\partial\Omega(X,Y,Z):=X\Omega(Y,Z)-\Omega([X,Y],Z)+\circlearrowleft
(X,Y,Z),
\end{gather}
where $  \circlearrowleft (X, Y,Z)$ denotes the summands obtained
from the two written ones by the cyclic permutation of the symbols
$X$, $Y$, $Z$.

In this paper, we are interested in the dif\/ferential ${\rm
\frak{sl}}(2)$-relative cohomology spaces
\[
{\mathrm H}^1_{\rm dif\/f}({\rm Vect_{Pol}}(\mathbb{R}),
\frak{sl}(2);{\cal
D}_{\lambda,\lambda+k})\qquad\text{and}\qquad{\mathrm H}^2_{\rm
dif\/f}({\rm Vect_{Pol}}(\mathbb{R}),{\rm \frak{sl}}(2);{\cal
D}_{\lambda,\lambda+k}).
\]
\begin{proposition}\label{invcoc} \null \quad {}

\begin{itemize}\itemsep=0pt
  \item [\rm i)] Any $1$-cocycle $c:{\rm Vect_{Pol}}(\mathbb{R})\rightarrow
  \mathcal{D}_{\lambda,\lambda+k}$ vanishing on ${\rm \frak{sl}}(2)$
coincides (up to a~scalar factor) with the transvectant
$J_{k+1}^{-1,\lambda}$ defined here and below by the formulae
\eqref{J}.
  \item [\rm ii)] Any  $2$-cocycle vanishing on ${\rm \frak{sl}}(2)$ is ${\rm
\frak{sl}}(2)$-invariant.
  \item [\rm iii)] Let $\Omega\in Z^2({\rm Vect_{Pol}}(\mathbb{R}),{\rm
\frak{sl}}(2);{\cal D}_{\lambda,\lambda+k})$. If $\Omega$ is a
${\rm \frak{sl}}(2)$-relative $2$-coboundary then (up to a~scalar
factor) we have $\Omega=\partial J_{k+1}^{-1,\lambda}$.
\end{itemize}
\end{proposition}

\begin{proof}
i) The 1-cocycle relation reads:
\[
c([X,Y])=L_X^{\mu}\circ c(Y)-c(Y)\circ
L_X^{\lambda}-L_Y^{\mu}\circ c(X)+c(X)\circ L_X^{\lambda}.
\]
Consider $X\in{\rm \frak{sl}}(2)$. Since $c(X)=0$, one easily sees
that
\begin{gather}\label{sltr1}
L_X^{\mu}c(Y)=c([X,Y])+c(Y)\circ L_X^{\lambda}.
\end{gather}
The equation (\ref{sltr1}) expresses  the
$\frak{sl}(2)$-invariance property of the bilinear map $c$. Thus,
according to Proposition~\ref{trans2}, the map $c$ coincides with
the transvectant $J_{k+1}^{-1,\lambda}$.

ii) Let $\Omega\in Z^2({\rm Vect_{Pol}}(\mathbb{R});{\cal
D}_{\lambda,\lambda+k})$. Since $ \Omega(X,Y)=0$ for all $X\in\rm
\frak{sl}(2) $, we deduce from the 2-cocycle relation~(\ref{2c})
that, for all $X\in{\rm \frak{sl}}(2)$ and $Y,\,Z\in{\rm
Vect_{Pol}}(\mathbb{R})$, we have
\begin{gather*}
X\Omega(Y,Z)-\Omega([X,Y],Z)-\Omega(Y,[X,Z])=0.
\end{gather*}
This last relation is nothing but the $\frak{sl}(2)$-invariance
property of the bilinear map $\Omega$.

iii) Let $\Omega=\partial b$. For all $X,\,Y\in{\rm
Vect_{Pol}}(\mathbb{R})$ we have
\[
\partial b(X,Y):= L_X^{\lambda,\lambda+k}
b(Y)-L_Y^{\lambda,\lambda+k}b(X)-b([X,Y]).
\]
Since $ \partial b(X,Y)=b(X)=0$ for all $X\in\frak {sl}(2) $ we
deduce that $b$ is ${\rm \frak{sl}}(2)$-invariant:
\[
L_X^{\lambda,\lambda+k} b(Y)=b([X,Y]).
\] According to Proposition \ref{trans2},
the space of $\mathfrak{sl}(2)$-invariant linear dif\/ferential
operator from \linebreak ${\rm Vect_{Pol}}(\mathbb{R})$ to ${\cal
D}_{\lambda,\lambda+k}$ vanishing on $\mathfrak{sl}(2)$ is one
dimensional and it is spanned by $J_{k+1}^{-1,\lambda}$. Thus, up
to a scalar factor, $b=J_{k+1}^{-1,\lambda}$.
Proposition~\ref{invcoc} is proved.
\end{proof}

The $\frak{sl}(2)$-trivial deformations are closely related to the
$\frak{sl}(2)$-relative cohomology  spaces ${\mathrm H}^1_{\rm
dif\/f}({\rm Vect_{Pol}}(\mathbb{R}), \frak{sl}(2);{\cal
D}_{\lambda,\lambda+k})$ and ${\mathrm H}^2_{\rm dif\/f}({\rm
Vect_{Pol}}(\mathbb{R}),{\rm \frak{sl}}(2);{\cal
D}_{\lambda,\lambda+k})$. Therefore, we will describe brief\/ly
these two spaces.

\subsection[The first cohomology space]{The f\/irst cohomology space}\label{sec3.1}

Note that, by Proposition \ref{invcoc}, we can describe the space
${\mathrm H}^1_{\rm dif\/f}({\rm Vect_{Pol}}(\mathbb{R}),{\rm
\frak{sl}}(2);{\cal D}_{\lambda,\lambda+k})$. This space is, in
fact, one-dimensional if and only if the corresponding
transvectant $J_{k+1}^{-1,\lambda}$ is a nontrivial ${\rm
\frak{sl}}(2)$-relative 1-cocycle, otherwise it is trivial.
However, this space  was computed  by Bouarroudj and Ovsienko, the
result is as follows:

\begin{theorem}[\cite{bo}]\label{thcoh1}
$\rm{dim}{\mathrm H}^1_{\rm dif\/f}({\rm
Vect_{Pol}}(\mathbb{R}),{\rm \frak{sl}}(2);{\cal
D}_{\lambda,\mu})=1$ if
\begin{gather*}
\mu-\lambda=2 \qquad \hbox{and} \qquad \lambda\neq-\tfrac{1}{2},\\
\mu-\lambda=3 \qquad \hbox{and} \qquad \lambda\neq-1,\\
\mu-\lambda=4 \qquad \hbox{and} \qquad \lambda\neq-\tfrac{3}{2},\\
\mu-\lambda=5 \qquad \hbox{and} \qquad \lambda=0,-4,\\
\mu-\lambda=6 \qquad \hbox{and} \qquad \lambda=-\tfrac{5\pm
\sqrt{19}}{2}.
\end{gather*}
Otherwise, ${\mathrm H}^1_{\rm dif\/f}({\rm
Vect_{Pol}}(\mathbb{R}),{\rm \frak{sl}}(2);{\cal
D}_{\lambda,\mu})=0$.
\end{theorem}

 These spaces ${\mathrm H}^1_{\rm dif\/f}({\rm
Vect_{Pol}}(\mathbb{R}),{\rm \frak{sl}}(2);{\cal
D}_{\lambda,\lambda+k})$ are generated by the cohomology classes
of the $\frak{sl}(2)$-relative 1-cocycles,
$C_{\lambda,\lambda+k}:{\rm
Vect_{Pol}}(\mathbb{R})\rightarrow\mathcal{D}_{\lambda,\lambda+k
}$ that are collected in the following~table.

\begin{table}[th]
\caption{}

\vspace{1mm}

\centerline{\begin{tabular}{|l|}\hline \tsep{1mm}\bsep{1mm}
$C_{\lambda,\lambda+2}(X,f)=X^{(3)}f$,\qquad
$\lambda\neq-\frac{1}{2}$\\\hline
$C_{\lambda,\lambda+3}(X,f)=X^{(3)}f'-\frac{\lambda}{2}X^{(4)}f$,\qquad $\lambda\neq-1$ \\
\hline \tsep{1mm}\bsep{1mm}
$C_{\lambda,\lambda+4}(X,f)=X^{(3)}f''-\frac{2\lambda+1}{2}X^{(4)}f'+\frac{\lambda(2\lambda+1)}{10}X^{(5)}f$,
\qquad $\lambda\neq-\frac{3}{2}$\\
\hline \tsep{1mm}\bsep{1mm} $C_{0,5}(X,f) =-3X^{(5)}f'+
15X^{(4)}f'' -10X^{(3)}f^{(3)}$\\ \hline \tsep{1mm}\bsep{1mm}
$C_{-4,1}(X,f)=28 X^{(6)}f+63X^{(5)}f'+ 45X^{(4)}f''
+10X^{(3)}f^{(3)}$\\ \hline \tsep{1mm}\bsep{1mm}
$C_{a_i,a_i+6}(X,f)=\alpha_iX^{(7)}f-14\beta_iX^{(6)}f'-126\gamma_iX^{(5)}f''
-210\tau_iX^{(4)}f^{(3)}+210X^{(3)}f^{4}$\\
\hline
\end{tabular}}\label{table1}

\vspace{2mm}

\noindent where
\begin{gather*}
a_1=-\tfrac{5+ \sqrt{19}}{2}, \quad \alpha_1=-\tfrac{22+
5\sqrt{19}}{4}, \quad \beta_1=\tfrac{31+ 7\sqrt{19}}{2}, \quad
\gamma_1=\tfrac{25+
7\sqrt{19}}{2}, \quad \tau_1=-2+\sqrt{19},\\
a_2=-\tfrac{5- \sqrt{19}}{2}, \quad \alpha_2=-\tfrac{22-
5\sqrt{19}}{4}, \quad  \beta_2=\tfrac{31- 7\sqrt{19}}{2},\quad
\gamma_2=\tfrac{25- 7\sqrt{19}}{2},\quad \tau_2=-2-\sqrt{19}.
\end{gather*}
\end{table}

For $X\frac{d}{dx}\in{\rm Vect_{Pol}}(\mathbb{R})$ and
$f{dx}^{\lambda}\in{\cal F}_\lambda$, we write
\begin{gather*}
C_{\lambda,\lambda+k
}(X\frac{d}{dx})(f{dx}^{\lambda})=C_{\lambda,\lambda+k
}(X,f){dx}^{\lambda+k}. \end{gather*}

The maps $C_{\lambda,\lambda+j}(X)$ are naturally extended  to
${\cal S}_\delta^n=\bigoplus_{j=0}^n{\cal F}_{\delta-j}$.

\subsection{The second cohomology space}\label{sec3.2}

Let $\mathfrak{g}$ be a Lie algebra, $\mathfrak{h}$ a subalgebra
of $\mathfrak{g}$ and $V$ a $\mathfrak{g}$-module, the {\it
cup-product} is def\/ined, for arbitrary linear maps $c_1 ,  c_2 :
\frak g\rightarrow \mathrm{End}(V)$, by
\begin{gather}
[\![c_1 , c_2]\!] : \ \ \frak g \otimes \frak g \rightarrow
\mathrm{End}(V),\qquad  {}[\![c_1 , c_2]\!] (x , y) = [c_1(x) ,
c_2(y)] + [c_2(x) , c_1(y)].\label{maurrer cartan1}
\end{gather}
It is easy to check that for any two
$\mathfrak{h}$-relative $1$-cocycles $c_1$ and $c_2 \in Z^1 (\frak
g ,\mathfrak{h} ; \mathrm{End}(V))$, the bilinear map $[\![c_1 ,
c_2]\!]$ is a $\mathfrak{h}$-relative $2$-cocycle. Moreover, if
one of the cocycles $c_1$ or $c_2$ is a  $\mathfrak{h}$-relative
1-coboundary, then $[\![c_1 , c_2]\!]$ is a
$\mathfrak{h}$-relative $2$-coboundary. Therefore, we naturally
deduce that the operation (\ref{maurrer cartan1}) def\/ines a
bilinear map
\begin{gather*}
\mathrm{H}^1 (\frak g ,\mathfrak{h};\mathrm{End}( V))\otimes
\mathrm{H}^1 (\frak g ,\mathfrak{h}; \mathrm{End}( V))\rightarrow
\mathrm{H}^2 (\frak g ,\mathfrak{h}; \mathrm{End}( V)).
\end{gather*}
Thus, by computing the cup-products of the 1-cocycles
$C_{\lambda,\lambda+k}$ generating the spaces
\[
{\mathrm H}^1_{\rm dif\/f}({\rm Vect_{Pol}}(\mathbb{R}),{\rm
\frak{sl}}(2);{\cal D}_{\lambda,\lambda+k}),
\] we can exhibit  explicit expressions of
some $\frak{sl}(2)$-relative 2-cocycles
\[
\Omega_{\lambda,\lambda+k}:{\rm
Vect_{Pol}}(\mathbb{R})\rightarrow{\cal D}_{\lambda,\lambda+k}.
\]

  For $X\frac{d}{dx},\,Y\frac{d}{dx}\in{\rm
Vect_{Pol}}(\mathbb{R})$ and $f{dx}^{\lambda}\in{\cal F}_\lambda$,
we write
\begin{gather*}
\Omega_{\lambda,\lambda+k
}\left(X\frac{d}{dx},Y\frac{d}{dx}\right)(f{dx}^{\lambda})=\Omega_{\lambda,\lambda+k
}(X,Y,f){dx}^{\lambda+k}. \end{gather*}
\begin{proposition}\label{prop2} {} \quad

\begin{itemize}\itemsep=0pt
  \item [\rm i)] The map $\Omega_{\lambda,\lambda+5}$ is defined by
\[
(\lambda+4)\Omega_{\lambda,\lambda+5}=2[\![C_{\lambda+2,\lambda+5},C_{\lambda,\lambda+2}]\!]
\]
is a nontrivial $\frak{sl}(2)$-relative $2$-cocycle if and only
if $\lambda\in\{0,-2,-4\}$. Moreover
\[
-2[\![C_{\lambda+3,\lambda+5},C_{\lambda,\lambda+3}]\!]=\lambda\Omega_{\lambda,\lambda+5}.
\]

  \item [\rm ii)] The map $\Omega_{\lambda,\lambda+6}$ is defined by
\[
(2\lambda+9)\Omega_{\lambda,\lambda+6}=-2[\![C_{\lambda+2,\lambda+6},C_{\lambda,\lambda+2}]\!]
\]
is a nontrivial $\frak{sl}(2)$-relative $2$-cocycle if and only
if $\lambda\in\left\{-\frac{5}{2}, -\frac{22- 5\sqrt{19}}{4},
-\frac{22+ 5\sqrt{19}}{4}\right\}$. Moreover,
\[
5(2\lambda+1)\Omega_{\lambda,\lambda+6}=-2(2\lambda+1)[\![C_{\lambda+3,\lambda+6},C_{\lambda,\lambda+3}]\!]=10
[\![C_{\lambda+4,\lambda+6},C_{\lambda,\lambda+4}]\!].
\]
\end{itemize}
\end{proposition}

\begin{proof}
By a straightforward computation we get
\begin{gather*}
\Omega_{\lambda,\lambda+5}(X,Y,f)=\big(X^{(4)}Y^{(3)}-X^{(3)}Y^{(4)}\big)f,\\
\Omega_{\lambda,\lambda+6}(X,Y,f)=\big(X^{(3)}Y^{(4)}-X^{(4)}Y^{(3)}\big)f'
-\tfrac{\lambda}{5}\big(X^{(3)}Y^{(5)}-X^{(5)}Y^{(3)}\big)f.
\end{gather*}
Moreover, we show also by a direct computation  that
\begin{gather*} 3\partial J_{6}^{-1,\lambda}=-
\lambda(\lambda^2+6\lambda+8)\Omega_{\lambda,\lambda+5},\qquad\text{and}\qquad
3\partial
J_{7}^{-1,\lambda}=(4\lambda^3+30\lambda^2+56\lambda+15)\Omega_{\lambda,\lambda+6}\!
\end{gather*} where
\begin{gather*}
J_6^{-1,\lambda}(X)(f) = 3 X^{(3)}f^{(3)}-\tfrac{9}{2}(\lambda+1)
X^{(4)}f'' +\tfrac{9}{10}(\lambda+1)(2\lambda+1)
X^{(5)}f'-\lambda\tfrac{2\lambda^2+3\lambda+1}{10}X^{(6)}f
\end{gather*}
and
\begin{gather*}
J_7^{-1,\lambda}(X)(f) =X^{(3)}f^{(4)}
-(2\lambda+3)X^{(4)}f^{(3)}+\tfrac{6\lambda^2+15\lambda+9}
{5}X^{(5)}f'' \\
\phantom{J_7^{-1,\lambda}(X)(f) =}{}
-\tfrac{4\lambda^3+12\lambda^2+11\lambda+3}{15}
X^{(6)}f'+\tfrac{\lambda(4\lambda^3+12\lambda^2+11\lambda+3)}{210}
X^{(7)}f.\end{gather*} Thus, we conclude by using
Proposition~\ref{invcoc}.
\end{proof}

\begin{proposition}\label{prop3}  The cup products
$[\![C_{\lambda+3,\lambda+7},C_{\lambda,\lambda+3}]\!]$ and
$[\![C_{\lambda+4,\lambda+7},C_{\lambda,\lambda+4}]\!]$ are
generically nontrivial $\frak{sl}(2)$-relative $2$-cocycles and they
are  cohomologous.
\end{proposition}
\begin{proof} The transvectant $J_8^{-1,\lambda}$ is given
by, for $X{d\over dx}\in{\rm Vect_{Pol}}(\mathbb{R})$ and
$fdx^\lambda\in\cal F_\lambda$,
\begin{gather*}
J_8^{-1,\lambda}(X)(f)= X^{(3)}f^{(5)}
-\tfrac{5}{2}(\lambda+2)X^{(4)}f^{(4)}+(\lambda+2)(2\lambda+3)X^{(5)}f^{(3)}\\
\phantom{J_8^{-1,\lambda}(X)(f)=}{} -\tfrac{1}{3}
(\lambda+1)(\lambda+2)(2\lambda+3)X^{(6)}f''+\tfrac{1}{42}(\lambda+1)(\lambda+2)(2\lambda+3)(2\lambda+1)X^{(7)}f'\!\!\\
\phantom{J_8^{-1,\lambda}(X)(f)=}{}-\tfrac{1}{840}\lambda(\lambda+1)(\lambda+2)(2\lambda+3)(2\lambda+1)X^{(8)}f.
\end{gather*}
Therefore, by a direct computation, we show that
\begin{gather*}\partial J_8^{-1,\lambda}(X,Y)(f) =\tfrac{1}{30}\lambda\Big(
(\lambda+1)(\lambda+2)(2\lambda+3)(2\lambda+11)+30\Big)X^{(3)}Y^{(6)}f\\
\qquad{}-\lambda(\lambda+2)\Big(-\tfrac{1}{60}(\lambda+1)(2\lambda+3)(2\lambda+1)+2\lambda+\tfrac{11}{2}\Big)
X^{(4)}Y^{(5)}f\\
 \qquad{}-\Big((\lambda+2)(2\lambda+3)\big(\tfrac{1}{3}(\lambda+1)(2\lambda+1)+
3\lambda+1\big)-5\lambda-1\Big)X^{(3)}Y^{(5)}f'\\
\qquad{}
+5\Big((\lambda+2)\big[\tfrac{1}{3}(\lambda+1)(2\lambda+3)+3\lambda+2\big]+2\lambda+1\Big)
 X^{(3)}Y^{(4)}f''-(X\leftrightarrow
Y).
\end{gather*}

Let us def\/ine $\Omega_{\lambda,\lambda+7}$ and
$\widetilde{\Omega}_{\lambda,\lambda+7}$ by
\[
\Omega_{\lambda,\lambda+7}=[\![C_{\lambda+3,\lambda+7},C_{\lambda,\lambda+3}]\!]\qquad{and}\qquad
\widetilde{\Omega}_{\lambda,\lambda+7}=[\![C_{\lambda+4,\lambda+7},C_{\lambda,\lambda+4}]\!].
\]
Thus,
\begin{gather*}\Omega_{\lambda,\lambda+7}(X,Y,f)= \left(\tfrac{-\lambda(2\lambda+7)(\lambda+8)}{20}X^{(5)}Y^{(4)}
-\tfrac{\lambda}{2}X^{(3)}Y^{(6)}\right)f+\tfrac{2\lambda^2+23\lambda+11}{10}X^{(5)}Y^{(3)}f'\\
\phantom{\Omega_{\lambda,\lambda+7}(X,Y,f)=}{}  +\tfrac{\lambda+11}{2}X^{(3)}Y^{(4)}f''-(X\leftrightarrow Y),\\
\widetilde{\Omega}_{\lambda,\lambda+7}(X,Y,f)=
\left(\tfrac{\lambda(2\lambda+1)}{10}X^{(3)}Y^{(6)}-
\tfrac{\lambda(\lambda+4)(2\lambda+1)}{20}X^{(4)}Y^{(5)}\right)f+
\tfrac{(\lambda-5)(2\lambda+1)}{10}X^{(3)}Y^{(5)}f'
\\
\phantom{\widetilde{\Omega}_{\lambda,\lambda+7}(X,Y,f)=}{}
+\tfrac{5-\lambda}{2}X^{(3)}Y^{(4)}f''-(X\leftrightarrow Y).
\end{gather*}
We exhibit some reals $a$, $b$ and $c$  such
that:
\begin{equation}\label{abc}
a\Omega_{\lambda,\lambda+7}+b\widetilde{\Omega}_{\lambda,\lambda+7}=c\partial J_8^{-1,\lambda}
.\end{equation} For $\lambda=-3$, we get
$\Omega_{\lambda,\lambda+7}=\widetilde{\Omega}_{\lambda,\lambda+7}$
and $\partial J_8^{-1,\lambda}=0$. For $\lambda=-6$, we get $a=1$
and $70c=5+11b$. For $\lambda\in\{-5,\,-{3\over2},\,-{1\over2}\}$,
we get
\[
a=0,\qquad b=1,\qquad\text{and}\qquad
210c=-8\lambda^3-60\lambda^2-70\lambda+45.
\] For
$\lambda\notin\{-6,\,-5,\,-3,\,-{3\over2},\,{1\over2}\}$, we get
\begin{gather*}
a=1,\qquad
b=\tfrac{4\lambda^3+48\lambda^2+161\lambda+117}{4\lambda^3+24\lambda^2+17\lambda-15}\qquad\text{and}\qquad
c=\tfrac{b(4\lambda^3+24\lambda^2+3\lambda-15)-4\lambda^3-48\lambda^2-147\lambda-33}{70(\lambda+3)}.
\end{gather*}
Now, it is easy to  show that, if
$4\lambda^3+48\lambda^2+161\lambda+117\neq0$ the 2-cocycle
$\Omega_{\lambda,\lambda+7}$ is nontrivial:
$\Omega_{\lambda,\lambda+7}\neq\partial J_8^{-1,\lambda}$.
Similarly, if $\lambda\notin\{-5,\,-{3\over2},\,{1\over2}\}$, we
show that $\widetilde{\Omega}_{\lambda,\lambda+7}$ is nontrivial.
\end{proof}

 \begin{proposition} The cup product
$[\![C_{\lambda,\lambda+4},C_{\lambda+4,\lambda+8}]\!]$ is a
nontrivial $\frak{sl}(2)$-relative $2$-cocycle.
\end{proposition}

\begin{proof} Let
$\Omega_{\lambda,\lambda+8}=[\![C_{\lambda,\lambda+4},C_{\lambda+4,\lambda+8}]\!].$
By a straightforward computation we show that
\begin{gather*}\Omega_{\lambda,\lambda+8}(X,Y,f)=
-\left(\tfrac{\lambda(2\lambda+1)(2\lambda+9)}{20}X^{(4)}Y^{(6)}
-\tfrac{\lambda(2\lambda+1)}{10}X^{(3)}Y^{(7)}\right)f\\
\phantom{\Omega_{\lambda,\lambda+8}(X,Y,f)=}{}
-\left(\tfrac{9(2\lambda+1)(2\lambda+9)}{20}X^{(5)}Y^{(4)}
-\tfrac{(2\lambda+1)(2\lambda-5)}{10}X^{(3)}Y^{(6)}\right)f'\\
\phantom{\Omega_{\lambda,\lambda+8}(X,Y,f)=}{}+\tfrac{18(1+\lambda)}{5}X^{(5)}Y^{(3)}f''-6X^{(4)}Y^{(3)}f^{(3)}-(X\leftrightarrow
Y).
\end{gather*}
As before, we show that this 2-cocycle
$\Omega_{\lambda,\lambda+8}$ is nontrivial:
$\Omega_{\lambda,\lambda+8}\neq\partial J_9^{-1,\lambda}$.
\end{proof}

Now, we collect in the following proposition some
$\frak{sl}(2)$-relative nontrivial 2-cocycles
$\Omega_{\lambda,\lambda+k}$ for $k=9,10$ and for singular values
of $\lambda$.

\begin{proposition} The following cup-products are $\frak{sl}(2)$-relative
nontrivial $2$-cocycles
\begin{alignat*}{3}
&\Omega_{0,9}=[\![C_{0,5},C_{5,9}]\!],&& \Omega_{-4,5}=[\![C_{-4,0},C_{0,5}]\!],&\\
& \Omega_{-8,1}=[\![C_{-8,-4},C_{-4,1}]\!],&&
\Omega_{a_i,a_i+9}=[\![C_{a_i,a_i+6},C_{a_i+6,a_i+9}]\!],&\\
&\Omega_{a_i-3,a_i+6}=[\![C_{a_i-3,a_i},C_{a_i,a_i+6}]\!],\qquad
&&
\Omega_{a_i,a_i+10}=[\![C_{a_i,a_i+6},C_{a_i+6,a_i+10}]\!],&\\
& \Omega_{a_i-4,a_i+6}=[\![C_{a_i-4,a_i},C_{a_i,a_i+6}]\!].&&&
\end{alignat*} \end{proposition}

\section{The general framework}\label{sec4}

In this section we def\/ine deformations of Lie algebra
homomorphisms and introduce the notion of miniversal deformations
over complete local algebras. Deformation theory of Lie algebra
homomorphisms was f\/irst considered with only one-parameter
deformation~\cite{ff2, nr2, r}. Recently, deformations of Lie
(super)algebras with multi-parameters were intensively studied
(see,  e.g.,~\cite{abbo, aalo2,  ro12, ro22}). Here we give an
outline of this theory.

\subsection[Infinitesimal deformations]{Inf\/initesimal deformations}\label{sec4.1}

Let $\rho_0:\frak g\to{\rm End}(V)$ be an action of a Lie algebra
$\frak g$ on a vector space $V$ and let $\frak h$ be a~subagebra
of $\frak g$. When studying $\frak h$-trivial deformations of the
$\frak g$-action $\rho_0$, one usually starts with inf\/initesimal
deformations
\[ \rho=\rho_0+tC,
\] where $C:\frak g\to{\rm
End}(V)$ is a linear map vanishing on $\frak h$ and $t$ is a
formal parameter. The homomorphism condition
\[
[\rho(x),\rho(y)]=\rho([x,y]),
\] where $x,y\in\frak g$, is
satisf\/ied in order 1 in $t$ if and only if $C$ is a $\frak
h$-relative 1-cocycle. That is, the map $C$ satisf\/ies
\begin{gather*}
[\rho_0(x),C(y)]-[\rho_0(y), C(x)]-C([x,y])=0.
\end{gather*}
Moreover, two $\frak h$-trivial inf\/initesimal deformations $
\rho=\rho_0+t\,C_1, $ and $ \rho=\rho_0+t\,C_2, $ are equivalents
if and only if $C_1-C_2$ is a $\frak h$-relative coboundary:
\begin{gather*}(C_1-C_2)(x)=[\rho_0(x),A]:=\partial
A(x),
\end{gather*}
where $A\in{\rm End}(V)^{\frak h}$ and $\partial$ stands for dif\/ferential
of cochains on $\frak g$ with values in $\mathrm{End}(V)$. So, the
space $\mathrm{H}^1(\frak g,\frak h;{\rm End}(V))$ determines and
classif\/ies the $\frak h$-trivial inf\/initesimal deformations up
to equivalence. (see,~e.g., \cite{Fuc2, nr2}). If
$\mathrm{H}^1(\frak g,\frak h;{\rm End}(V))$ is multi-dimensional,
it is natural to consider multi-parameter $\frak h$-trivial
deformations. More precisely, if $\mathrm{dim}{\mathrm H}^1(\frak
g,\frak h;{\rm End}(V))=m$, then choose $\frak h$-relative
1-cocycles $C_1,\ldots,C_m$ representing a basis of ${\mathrm
H}^1(\frak g,\frak h;{\rm End}(V))$ and consider the $\frak
h$-trivial inf\/initesimal deformation
\begin{gather*}
\rho=\rho_0+\sum_{i=1}^m t_i  C_i,
\end{gather*}
with independent parameters $t_1,\ldots,t_m$.

In our study, we are interested in the inf\/initesimal ${\rm
\frak{sl}}(2)$-trivial deformation of the ${\rm
Vect_{Pol}}(\mathbb{R})$-action on ${\cal
S}_\delta^n=\bigoplus_{j=0}^n{\cal F}_{\delta-j}$, the space of
symbols of dif\/ferential operators, where $n\in\mathbb{N}$ and
$\delta\in\mathbb{R}$. Thus, we  consider the $\frak
{sl}(2)$-relative cohomology space ${\mathrm H}^1_{\rm
dif\/f}({\rm Vect_{Pol}}(\mathbb{R}),\frak{ sl}(2);{\cal D})$. Any
inf\/initesimal $\frak{sl}(2)$-trivial deformation is then of the
form
\begin{gather}\label{InDefGen2}{\cal L}_X=L_X+{\cal
L}^{(1)}_X, \end{gather} where $L_X$ is the Lie derivative of
${\cal S}_\delta^n$ along the vector f\/ield $X\frac{d}{dx}$
def\/ined by (\ref{Lieder2}), and
\begin{gather}
\label{l1} {\cal L}_X^{(1)}=\sum_\lambda\sum_{j=2}^6{}t_{\lambda,
\lambda+j}  C_{\lambda,\lambda+j}(X) \end{gather}
 and where $t_{\lambda,\lambda+j}$  are independent parameters,
$\delta-\lambda\in\mathbb{N}$,
$\delta-n\leq\lambda,\lambda+j\leq\delta$ and the $\frak{
sl}(2)$-relative 1-cocycles $C_{\lambda,\lambda+j}$ are def\/ined
in Table~\ref{table1}.

Note that for $(j,\lambda)=(2, -\frac{1}{2}), (3, -1), (4,
-\frac{3}{2})$, or $j=5$ and $\lambda\notin\{0,-4\}$ or $j=6$ and
$\lambda\neq -\frac{5\pm \sqrt{19}}{2}$ we have $
C_{\lambda,\lambda+j}=0$, then there are no corresponding
parameters $t_{\lambda,\lambda+j}$.

\subsection{Integrability conditions}\label{sec4.2}

Consider the problem of integrability of inf\/initesimal
deformations. Starting with the in\-f\/i\-ni\-te\-si\-mal
deformation (\ref{InDefGen2}), we look for a formal series
\begin{gather}
\label{formal1} {\cal L}_X=L_X+{\cal L}^{(1)}_X+{\cal
L}^{(2)}_X+{\cal L}^{(3)}_X+\cdots,
\end{gather}
where ${\cal L}_X^{(k)}$ is an homogenous polynomial of degree $k$
in the parameters $(t_{\lambda,\lambda+j})$ and with
coef\/f\/icients in ${\cal D}$ such that ${\cal L}_X^{(k)}=0$ if
$X\frac{d}{dx}\in{\rm \frak{sl}}(2)$. This formal series
(\ref{formal1}) must satisfy the homomorphism condition in any
order in the parameters $(t_{\lambda,\lambda+j})$ \begin{gather}
\label{Homom} [{\cal L}_X,{\cal L}_Y]={\cal L}_{[X,Y]}.
\end{gather}
The homomorphism condition (\ref{Homom}) gives the following
(Maurer--Cartan) equations
\begin{gather}\label{MCG}\partial{\cal L}^{(k)}=-\tfrac{1}{2} \sum_{i+j=k}[\![{\cal
L}^{(i)},{\cal L}^{(j)}]\!]. \end{gather}

However, quite often the above problem has no solution. Note here
that the right side of~(\ref{MCG}) must be a coboundary of a
1-cochain vanishing on ${\rm \frak{sl}}(2)$, so, the obstructions
for integrability of inf\/initesimal deformations belong to the
second ${\rm \frak{sl}}(2)$-relative cohomology space ${\mathrm
H}^2_{\rm dif\/f}({\rm Vect_{Pol}}(\mathbb{R}),{\rm
\frak{sl}}(2);{\cal D})$.

Following \cite{ff2} and \cite{aalo2}, we will impose extra
algebraic relations on the parameters $(t_{\lambda,\lambda+j})$.
Let~${\cal R}$ be an ideal in
$\mathbb{C}[[t_{\lambda,\lambda+j}]]$ generated by some set of
relations, the quotient
\begin{gather}
\label{TrivAlg2} {\cal
A}=\mathbb{C}[[t_{\lambda,\lambda+j}]]/{\cal R}
\end{gather}
is a complete local algebra with unity, and one can speak about
deformations with base~${\cal A}$, see~\cite{ff2} for details.

Given an inf\/initesimal deformation (\ref{InDefGen2}), one can
always consider it as a deformation with base (\ref{TrivAlg2}),
where $\cal R$ is the ideal  generated by all the quadratic
monomials. Our aim is to f\/ind $\cal A$ which is big as possible,
or, equivalently, we look for relations on the parameters
$(t_{\lambda,\lambda+j})$ which are necessary and suf\/f\/icient
for integrability (cf.~\cite{abbo,aalo2}).

\subsection{Equivalence and the miniversal deformation}\label{sec4.3}

The notion of equivalence of deformations over complete local
algebras has been considered in~\cite{f2}.

\begin{definition}  Two deformations, $\rho$ and $\rho'$ with the same
base $\cal A$ are called equivalent if there exists an inner
automorphism $\Psi$ of the associative algebra ${\rm
End}(V)\otimes\cal A$ such that
\[
\Psi\circ\rho=\rho'\hbox{ and }\Psi(\mathbb{I})=\mathbb{I},
\]
where $\mathbb{I}$ is the unity of the algebra ${\rm
End(V)}\otimes\cal A$.
\end{definition}

The following notion of miniversal deformation is fundamental. It
assigns to a~$\frak g$-module $V$ a~canonical commutative
associative algebra $\cal A$ and a canonical deformation with base
$\cal A$.

\begin{definition} A deformation $\rho$ with base $\cal A$ is called
miniversal, if
\begin{itemize}\itemsep=0pt
  \item [(i)] for any other deformation, $\rho'$ with base (local)
$\cal A'$, there exists a  homomorphism $\psi:{\cal A}'\to{\cal
A}$ satisfying $\psi(1)=1$, such that
\[
\rho=(\mathbb{I}d\otimes\psi)\circ\rho'.
\]
  \item [(ii)] in the notations of (i), if $\mathcal{A}$ is inf\/initesimal
  then  $\psi$ is unique.
\end{itemize}
If $\rho$ satisf\/ies only the condition (i), then it is called
versal.
\end{definition}
 We refer to~\cite{ff2} for a construction of miniversal
deformations of Lie algebras and to~\cite{aalo2} for miniversal
deformations of $\frak g$-modules.

\section{Second-order integrability conditions}\label{sec5}

In this section we obtain the integrability conditions for the
inf\/initesimal deformation~(\ref{InDefGen2}).  Assume that the inf\/initesimal deformation (\ref{InDefGen2})
can be integrated to a formal deformation
\[
 {\cal L}_X=L_X+{\cal
L}^{(1)}_X+{\cal L}^{(2)}_X+{\cal L}^{(3)}_X+\cdots,
\]
 where ${\cal L}^{(1)}_X$ is given
by (\ref{l1}) and ${\cal L}^{(2)}_X$ is a quadratic polynomial in
$t$ whose coef\/f\/icients are elements of ${\cal D}$ vanishing on
$\frak{ sl}(2)$. We compute the conditions for the second-order
terms~${\cal L}^{(2) }$. The homomorphism condition
\[
[{\cal L}_X,{\cal L}_Y]={\cal L}_{[X,Y]},
\]
gives for the second-order terms the following (Maurer--Cartan)
equation
\begin{gather}\label{MC1} \partial {\cal L}^{(2)}=-\tfrac{1}{2} [\![{\cal
L}^{(1)},{\cal L}^{(1)}]\!].
\end{gather} The right hand side of
(\ref{MC1}) is a cup-product of ${\rm \frak{sl}}(2)$-relative
1-cocycles, so it is automatically a~${\rm \frak{sl}}(2)$-relative
2-cocycle. More precisely, the equation (\ref{MC1}) can be
expressed as follows
\begin{gather}\label{cap}\partial{\cal
L}^{(2)}=-\tfrac{1}{2}[\![\sum_\lambda\sum_{j=2}^6{}t_{\lambda,
\lambda+j}\,C_{\lambda,\lambda+j},\sum_\lambda\sum_{j=2}^6{}t_{\lambda,\lambda+j}
\,C_{\lambda,\lambda+j}]\!],\end{gather} therefore, let us
consider the ${\rm \frak{sl}}(2)$-relative 2-cocycles
$B_{\lambda,\lambda+k}\in Z^2_{\rm dif\/f}({\rm
Vect_{Pol}}(\mathbb{R}),{\rm \frak{sl}}(2),{\cal
D}_{\lambda,\lambda+k})$, for $k=4,\dots,10,$ def\/ined by
\[
B_{\lambda,\lambda+k}=-\sum_{j=2}^k{}t_{\lambda+j,\lambda+k}{}t_{\lambda,\lambda+j}
\,[\![C_{\lambda+j,\lambda+k},C_{\lambda,\lambda+j}]\!].
\]
It is easy to see that $B_{\lambda,\lambda+4}=0 $. The second
order integrability conditions are determined by the fact that any map
2-cocycles $B_{\lambda,\lambda+k}$, for $k = 5, \dots, 10$, must
be a ${\rm \frak{sl}}(2)$-relative 2-coboundary. More precisely,
$B_{\lambda,\lambda+k}$ must coincide, up to a scalar factor, with
$\partial J_{k+1}^{-1,\lambda}$. We split these conditions into two
family which we explain in the two following propositions. Let us
f\/irst consider the following functions in $t$ where $t$ is the
family of parameters $(t_{\lambda,\lambda+j})$
\begin{gather*} \omega_{\lambda,\lambda+5}(t)=-\tfrac{\lambda+4}{2}
t_{\lambda,\lambda+2}t_{\lambda+2,\lambda+5}+
\tfrac{\lambda}{2}t_{\lambda,\lambda+3}t_{\lambda+3,\lambda+5},\\
\omega_{\lambda,\lambda+6}(t) =\tfrac{2\lambda+9}{2}
t_{\lambda,\lambda+2}t_{\lambda+2,\lambda+6}+\tfrac{5}{2}t_{\lambda,\lambda+3}t_{\lambda+3,\lambda+6}
-\tfrac{2\lambda+1}{2}
t_{\lambda,\lambda+4}t_{\lambda+4,\lambda+6},\\
\omega_{\lambda,\lambda+7}(t) =\left\{\begin{array}{ll} -{c\over
b}\,t_{\lambda,\lambda+4}t_{\lambda+4,\lambda+7},\quad\text{if}\quad
b\neq0\\
-c\,t_{\lambda,\lambda+3}t_{\lambda+3,\lambda+7},\quad\text{if}\quad
b=0\end{array}\right.,
\quad\text{if}\quad\lambda\notin\{-6,\,0\},\\
\omega_{0,7}(t)=-\tfrac{1}{7}\big(\tfrac{11}{10}t_{0,3}t_{3,7}+\tfrac{1}{2}t_{0,4}t_{4,7}+3t_{0,5}t_{5,7}\big),\\[4pt]
\omega_{-6,1}(t)=\tfrac{1}{14}\big(t_{-6,-3}t_{-3,1}-6t_{-6,-4}t_{-4,1}+\tfrac{11}{5}t_{-6,-2}t_{-2,1}\big)
,\\[4pt]
\omega_{0,8}(t)=\tfrac{2}{11}t_{0,5}t_{5,8},\\
\omega_{-7,1}(t)=\tfrac{2}{15}t_{-7,-4}t_{-4,1}.
\end{gather*}
These functions $\omega_{\lambda,\lambda+k}(t)$, $k=5,6,7$, will
appear as coef\/f\/icients for some maps  from ${\cal F}_\lambda$
to ${\cal F}_{\lambda+k}$ and they will be used in the expressions
of integrability conditions. More precisely, we will show that the
second term ${\cal L}^{(2)}$ is of the form ${\cal
L}^{(2)}=\sum_{\lambda,k}\omega_{\lambda,\lambda+k}(t)J_{k+1}^{-1,\lambda}$.

\begin{proposition}\label{pr3} For $k=5,6,7$, we have the following
second-order integrability  conditions of the infinitesimal
deformation~\eqref{InDefGen2}
\begin{alignat}{3}
&\omega_{\lambda,\lambda+5}(t)=0\qquad &&\text{if} \ \ \lambda\in\{0,-2,-4\},& \nonumber\\
& \omega_{\lambda,\lambda+6}(t) =0\qquad &&\text{if} \ \ \lambda\in\left\{-\tfrac{5\pm \sqrt{19}}{2},-\tfrac{5}{2}\right\},&\nonumber\\
& b\,t_{\lambda,\lambda+3}t_{\lambda+3,\lambda+7}-a\,t_{\lambda,\lambda+4}t_{\lambda+4,\lambda+7}
= 0\qquad &&\text{if} \ \ \lambda\notin\{0,-2,-4,-6\},&\nonumber\\
& 10 t_{-2,0}t_{0,5}-t_{-2,1}t_{1,5}-\tfrac{1}{3}t_{-2,2}t_{2,5}=0,&&&\nonumber\\
& 10 t_{-4,1}t_{1,3}+t_{-4,-1}t_{-1,3}+3\,t_{-4,0}t_{0,3}=0.&&&\label{k567}
\end{alignat}
where $a$ and $b$ are defined by \eqref{abc}.
\end{proposition}

\begin{proof} { 1)} For $k=5$, we have
\[
B_{\lambda,\lambda+5}=-t_{\lambda,\lambda+2}t_{\lambda+2,\lambda+5}
[\![C_{\lambda+2,\lambda+5},C_{\lambda,\lambda+2}]\!]-t_{\lambda,\lambda+3}t_{\lambda+3,\lambda+5}
[\![C_{\lambda+3,\lambda+5},C_{\lambda,\lambda+5}]\!],
\]
hence, according to Proposition \ref{prop2}, we have
\[
B_{\lambda,\lambda+5}=\omega_{\lambda,\lambda+5}(t)\Omega_{\lambda,\lambda+5}.
\]
Thus, by Proposition~\ref{prop2},  the ${\rm
\frak{sl}}(2)$-relative 2-cocycle $\Omega_{\lambda,\lambda+5}$ is
nontrivial if and only if $\lambda\in\{0,-2,-4\}$. Hence, for
$\lambda\in\{0,-2,-4\}$, the condition $
\omega_{\lambda,\lambda+5}(t)=0 $ holds.

2)  For $k=6$, as before, we have{\samepage
\[
B_{\lambda,\lambda+6}=\omega_{\lambda,\lambda+6}(t)\Omega_{\lambda,\lambda+6}.
\]
Thus, if $\lambda\in\left\{-\frac{5\pm
\sqrt{19}}{2},-\frac{5}{2}\right\}$ the condition $
\omega_{\lambda,\lambda+6}(t)=0 $ must be satisf\/ied.}

3) Let $k=7$. Note that, hereafter, some singular values of the
parameter $\lambda$ appear because  the ${\rm
\frak{sl}}(2)$-relative 2-cocycles $C_{\lambda,\lambda+5}$ exist
only for $\lambda=0,-4$ and $C_{\lambda,\lambda+6}$ exist only for
$\lambda=-\frac{5\pm \sqrt{19}}{2}$.

i) For $\lambda\notin\{0,-2,-4,-6\}$, we have
\[
B_{\lambda,\lambda+7}=-t_{\lambda,\lambda+3}t_{\lambda+3,\lambda+7}\Omega_{\lambda,\lambda+7}
-t_{\lambda,\lambda+4}t_{\lambda+4,\lambda+7}\widetilde{\Omega}_{\lambda,\lambda+7}.
\]
Therefore, the following conditions follow from
Proposition~\ref{prop3}
\[
b\,t_{\lambda,\lambda+3}t_{\lambda+3,\lambda+7}-a\,t_{\lambda,\lambda+4}t_{\lambda+4,\lambda+7}=0.\]
Indeed, according to the equation (\ref{abc}), if $b\neq0$, we have
\[
B_{\lambda,\lambda+7}=-\tfrac{1}{b}(b t_{\lambda+3,\lambda+7}-a t_{\lambda,\lambda+4}t_{\lambda+4,\lambda+7})
\Omega_{\lambda,\lambda+7}+\omega_{\lambda,\lambda+7}(t)\partial
J_8^{-1,\lambda}
\]
 and if $b=0$, we have
 \[
 B_{\lambda,\lambda+7}=- t_{\lambda,\lambda+4}t_{\lambda+4,\lambda+7}
\widetilde{\Omega}_{\lambda,\lambda+7}+\omega_{\lambda,\lambda+7}(t)\partial
J_8^{-1,\lambda}.
\]

ii) By a direct computation, we show that
\[
B_{0,7}=\omega_{0,7}(t)\partial J_8^{-1,0}\qquad\text{and}\qquad
B_{-6,1}=\omega_{-6,1}(t)\partial J_8^{-1,-6}.
\]
Hence, there are no conditions on $B_{0,7}$ and $B_{-6,1}$.

iii) For the other two singular values of $\lambda$ we have
\begin{gather*}
B_{-2,5} =
-t_{-2,0}t_{0,5}[\![C_{0,5},C_{-2,0}]\!]-t_{-2,1}t_{1,5}
\Omega_{-2,5}-t_{-2,2}t_{2,5}\widetilde{\Omega}_{-2,5},\\
B_{-4,3} =
-t_{-4,1}t_{1,3}[\![C_{1,3},C_{-4,1}]\!]-t_{-4,-1}t_{-1,3}
\Omega_{-4,3}-t_{-4,0}t_{0,3}\widetilde{\Omega}_{-4,3}.
\end{gather*}

More precisely, we get
\begin{gather*}
B_{-2,5} = \big(10 t_{-2,0}t_{0,5}-t_{-2,1}t_{1,5}-\tfrac{1}{3}t_{-2,2}t_{2,5}\big)\Omega_{-2,5}+\omega_{-2,5}(t)\partial
J_8^{-1,-2},\\
B_{-4,3} = -(10 t_{-4,1}t_{1,3}+t_{-4,-1}t_{-1,3}+3 t_{-4,0}t_{0,3})\Omega_{-4,3}+\omega_{-4,3}(t)\partial
J_8^{-1,-4}.
\end{gather*}
So, the following integrability conditions become again from
Proposition~\ref{prop3}:
\begin{gather*}
10\,t_{-2,0}t_{0,5}-t_{-2,1}t_{1,5}-\tfrac{1}{}
t_{-2,2}t_{2,5}=10 t_{-4,1}t_{1,3}+t_{-4,-1}t_{-1,3}+3 t_{-4,0}t_{0,3}=0.\tag*{\qed}
\end{gather*}\renewcommand{\qed}{}
\end{proof}

\begin{proposition}\label{pr4}For $k=8,9,10$, we have the following
second-order integrability  conditions of the infinitesimal
deformation~\eqref{InDefGen2}, where in the first line
$\lambda\notin\{0,-3,-4,-7\}$,
\begin{gather}
t_{\lambda,\lambda+4}t_{\lambda+4,\lambda+8}=
11t_{0,4}t_{4,8}+10t_{0,5}t_{5,8}=0,\nonumber\\
t_{-3,1}t_{1,5}-10t_{-3,0}t_{0,5}=
t_{-4,0}t_{0,4}+10t_{-4,1}t_{1,4}=0,\nonumber\\
11t_{-7,-3}t_{-3,1}-10t_{-7,-4}t_{-4,1}=
t_{a_i,a_i+6}t_{a_i+6,a_i+8}=0,\nonumber\\
t_{a_i-2,a_i}t_{a_i,a_i+6}= t_{0,5}t_{5,9}=0,\nonumber\\
t_{-4,0}t_{0,5}-t_{-4,1}t_{1,5}= t_{-8,-4}t_{-4,1} = 0,\nonumber\\
t_{a_i,a_i+6}t_{a_i+6,a_i+9}= t_{a_i-3,a_i}t_{a_i,a_i+6}=0,\nonumber\\
t_{a_i,a_i+6}t_{a_i+6,a_i+10}=
t_{a_i-4,a_i}t_{a_i,a_i+6}=0.\label{k89}
\end{gather}
\end{proposition}

\begin{proof}
 1) For $k=8$, we f\/irst recall that the cup-product
$\Omega_{\lambda,\lambda+8}=[\![C_{\lambda,\lambda+4},C_{\lambda+4,\lambda+8}]\!]$
is a~${\rm \frak{sl}}(2)$-relative nontrivial 2-cocycle. Moreover,
for $\lambda\notin\{0,-3,-4,-7,a_1,a_2,a_1-2,a_2-2\}$, we have
\[
B_{\lambda,\lambda+8}=-
t_{\lambda,\lambda+4}t_{\lambda+4,\lambda+8}\Omega_{\lambda,\lambda+8}.
\]

For the  singular values, we easily check  that
\begin{alignat*}{3}
&[\![C_{5,8},C_{0,5}]\!]=\tfrac{10}{11}\Omega_{0,8}+\tfrac{2}{11}
\partial J_9^{-1,0},\qquad &&[\![C_{0,5},C_{-3,0}]\!]=10\Omega_{-3,5},&\\
&[\![C_{-4,1},
C_{-7,-4}]\!]=\tfrac{10}{11}\Omega_{-7,1}+\tfrac{2}{15}\partial
J_9^{-1,-7},\qquad &&[\![C_{1,4},C_{-4,1}]\!]=-10\Omega_{-4,4}&
\end{alignat*}
and we show that $[\![C_{a_i+6,a_i+8},C_{a_i,a_i+6}]\!]$, and
$[\![C_{a_i,a_i+6},C_{a_i-2,a_i}]\!]$ are  also  nontrivial
2-cocycles.

Thus, we deduce all integrability conditions corresponding to the
case $k=8$.

2) For $k=9$, the integrability conditions follow from the fact
that any corresponding cup-product of 1-cocycle is nontrivial.
Moreover, we have only singular cases and we also show that
$[\![C_{1,5},C_{-4,1}]\!]=-[\![C_{0,5},C_{-4,0}]\!]$.

3) For $k=10$ and $\lambda\neq a_i, a_i-4$ we have
$B_{\lambda,\lambda+10}=0.$ For $\lambda= a_i, a_i-4$ we have
\begin{gather*}
B_{a_i,a_i+10}=-t_{a_i,a_i+6}t_{a_i+6,a_i+10}[\![C_{a_i,a_i+6},C_{a_i+6,a_i+10}]\!]
=-t_{a_i,a_i+6}t_{a_i+6,a_i+10}\Omega_{a_i,a_i+10},\\
B_{a_i-4,a_i+6}=-t_{a_i-4,a_i}t_{a_i,a_i+6}[\![C_{a_i-4,a_i},C_{a_i,a_i+6}]\!]
=-t_{a_i-4,a_i}t_{a_i,a_i+6}\Omega_{a_i-4,a_i+6}.
\end{gather*}
Like in the previous case we prove that the 2-cocycles
$\Omega_{a_i,a_i+10}$ and $\Omega_{a_i-4,a_i+6}$ are nontrivial
and then we deduce the corresponding integrability conditions.
\end{proof}

Our main result in this section is the following

\begin{theorem}
\label{Main1}  The conditions \eqref{k567} and \eqref{k89} are
necessary and sufficient for second-order integrability of the
$\frak{sl}(2)$-trivial infinitesimal
deformation~\eqref{InDefGen2}.
\end{theorem}

\begin{proof} Of course, these conditions are necessary as, it was shown
in Proposition~\ref{pr3} and Proposition~\ref{pr4}. Now, under
these conditions, the second term $\mathcal{L}^{(2)}$ of the the
$\frak{sl}(2)$-trivial inf\/initesimal
deformation~(\ref{InDefGen2}) is a solution of the Maurer--Cartan
equation~(\ref{cap}). This solution is def\/ined up to a
1-coboundary and it has been shown in~\cite{ff2,aalo2} that
dif\/ferent choices of solutions of the Maurer--Cartan equation
correspond to equivalent deformations. Thus, we can always choose
\begin{gather*}
{\cal L}^{(2)}= \tfrac{1}{2}\sum_{\lambda\neq0,-2,-4}
\omega_{\lambda,\lambda+5}(t)J_{6}^{-1,\lambda}+
\tfrac{1}{2}\sum_{\lambda\neq
a_i,-\tfrac{5}{2}}\omega_{\lambda,\lambda+6}(t)J_{7}^{-1,\lambda}
\nonumber\\
\phantom{{\cal L}^{(2)}=}{} + \tfrac{1}{2}\sum_{\lambda}
\omega_{\lambda,\lambda+7}(t)J_{8}^{-1,\lambda}+
\tfrac{1}{2}\sum_{\lambda=0,-7}\omega_{\lambda,\lambda+8}(t)J_{9}^{-1,\lambda}.
\end{gather*}
Of course, any $t_{\lambda,\lambda+k}$ appears in the expressions
of ${\cal L}^{(1)}$ or ${\cal L}^{(2)}$ if and only if
$\delta-\lambda$ and $k$ are integers satisfying
$\delta-n\leq\lambda,\lambda+k\leq\delta$. Theorem~\ref{Main1} is
proved.
\end{proof}

\section{Third and fourth-order integrability conditions}\label{sec6}

\subsection[Computing the third-order Maurer-Cartan equation]{Computing the third-order Maurer--Cartan equation}\label{sec6.1}

Now we reconsider the formal deformation (\ref{formal1}) which is
a formal power series in the parame\-ters~$t_{\lambda,\lambda+j}$
with coef\/f\/icients in ${\cal D}$. We suppose that the
second-order integrability conditions are satisf\/ied. So, the
third-order terms of (\ref{formal1}) are solutions of the
(Maurer--Cartan) equation
\begin{gather}
\label{MC2} \partial {\cal L}^{(3)}=-\tfrac{1}{2}
\sum_{i+j=3}[\![{\cal L}^{(i)},{\cal L}^{(j)}]\!].
\end{gather}
As in the previous section we can write
\begin{gather}\label{MC4}\partial{\cal L}^{(3)}=-\tfrac{1}{2} \sum_{k,\lambda} E_{\lambda,\lambda+k},
\end{gather} where $E_{\lambda,\lambda+k}$ are maps from
${\rm Vect_{Pol}}(\mathbb{R})\times{\rm Vect_{Pol}}(\mathbb{R})$
to ${\cal D}_{\lambda,\lambda+k}$.  The third-order term
$\mathcal{L}^{(3)}$ of the $\frak{sl}(2)$-trivial formal
deformation (\ref{formal1}) is a solution of (\ref{MC4}). So, the
2-cochains $E_{\lambda,k}$ must satisfy $E_{\lambda,k}=\partial
J_{k+1}^{-1,\lambda}$ and then the third-order integrability
conditions are deduced from this fact.

It is easy to see that $E_{\lambda,\lambda+k}=0 $ for $k\leq6$ or
$k\geq13$, so we compute successively the $E_{\lambda,\lambda+k}$
for $k=7,\dots,12 $ and we resolve $E_{\lambda,\lambda+k}=\partial
J_{k+1}^{-1,\lambda}$ to get the corresponding third-order
integrability conditions.

Here, we mention  that the maps $E_{\lambda,\lambda+k}$ are
2-cochains, but they are not necessarily 2-cocycles because they
are not cup-products of 1-cocycles like the maps
$B_{\lambda,\lambda+k}$. Indeed, ${\cal L}^{(2)}$ is not
necessarily a 1-cocycle.

\subsection{Third-order integrability conditions}\label{sec6.2}

\begin{proposition}\label{pr5}For $k=7,8$, we have the following
third-order integrability  conditions of the infinitesimal
deformation~\eqref{InDefGen2}, for all $\lambda$
\begin{gather}
t_{\lambda,\lambda+2} \omega_{\lambda+2,\lambda+7}(t)=
\omega_{\lambda,\lambda+5}(t) t_{\lambda+5,\lambda+7}=0,\nonumber\\
t_{\lambda,\lambda+2}\,\omega_{\lambda+2,\lambda+8}(t)=
\omega_{\lambda,\lambda+6}(t) t_{\lambda+6,\lambda+8}=0,\nonumber\\
t_{\lambda,\lambda+3}\omega_{\lambda+3,\lambda+8}(t)=
\omega_{\lambda,\lambda+5}(t)
t_{\lambda+5,\lambda+8}=0.\label{thirdk7}
\end{gather}
\end{proposition}

\begin{proof}
 For $k=7$ and $\lambda\notin\{0,-2,-4,-6\}$ we have
\begin{gather*}E_{\lambda,\lambda+7}=t_{\lambda,\lambda+2}
\omega_{\lambda+2,\lambda+7}(t)
[\![J_{6}^{-1,\lambda+2},\,C_{\lambda,\lambda+2}]\!]+t_{\lambda+5,\lambda+7}
\omega_{\lambda,\lambda+5}(t)[\![C_{\lambda+5,\lambda+7},\,J_{6}^{-1,\lambda}]\!],
\\
E_{0,7}=t_{0,2} \omega_{2,7}(t)
[\![J_{6}^{-1,2},\,C_{0,2}]\!],\qquad E_{-6,1}=t_{-1,1}
\omega_{-6,-1}(t)[\![C_{-1,1},\,J_{6}^{-1,-6}]\!]
\end{gather*}
and
\[
E_{-2,5}=E_{-4,3}=0.
\]

By a direct computation, we show that the three maps
$[\![J_{6}^{-1,\lambda+2}, C_{\lambda,\lambda+2}]\!]$,
$[\![C_{\lambda+5,\lambda+7}, J_{6}^{-1,\lambda}]\!]$ and
$\partial J_8^{-1,\lambda}$ are linearly independent, for all
$\lambda$. Thus,
\begin{gather*}
t_{\lambda,\lambda+2}
\omega_{\lambda+2,\lambda+7}(t)=0 \qquad\text{for}\quad\lambda\neq-2,-4,-6,\\
t_{\lambda+5,\lambda+7}
\omega_{\lambda,\lambda+5}(t)=0 \qquad\text{for}\quad\lambda\neq0,-2,-4.
\end{gather*}
But, under the second-order integrability conditions:
$\omega_{\lambda,\lambda+5}(t)=0$ for $\lambda\in\{0,-2,-4\}$, the
conditions
\[t_{\lambda,\lambda+2}
\omega_{\lambda+2,\lambda+7}(t)=\omega_{\lambda,\lambda+5}(t)t_{\lambda+5,\lambda+7}=0
\]
hold for all $\lambda$.

 Now, for $k=8$ and $\lambda\notin\{ a_1,
a_2,-\frac{5}{2},a_1-2,a_2-2,-\frac{9}{2},0,-2,-4,-3,-5,-7\}$ we
have
\begin{gather*}E_{\lambda,\lambda+8} =t_{\lambda,\lambda+2}
\omega_{\lambda+2,\lambda+8}(t) [\![J_{7}^{-1,\lambda+2},
C_{\lambda,\lambda+2}]\!]+t_{\lambda+6,\lambda+8}
\omega_{\lambda,\lambda+6}(t)[\![C_{\lambda+6,\lambda+8}, J_{7}^{-1,\lambda}]\!]\\
\phantom{E_{\lambda,\lambda+8} =}{} +t_{\lambda,\lambda+3}
\omega_{\lambda+3,\lambda+8}(t) [\![J_{6}^{-1,\lambda+3},
C_{\lambda,\lambda+3}]\!]+t_{\lambda+5,\lambda+8}
\omega_{\lambda,\lambda+5}(t)[\![C_{\lambda+5,\lambda+8},
J_{6}^{-1,\lambda}]\!].
\end{gather*}
As before, we show that
\begin{gather*}
t_{\lambda,\lambda+2}
\omega_{\lambda+2,\lambda+8}(t)=t_{\lambda,\lambda+3}
\omega_{\lambda+3,\lambda+8}(t)=t_{\lambda+5,\lambda+8}
\omega_{\lambda,\lambda+5}(t)=t_{\lambda+6,\lambda+8}
\omega_{\lambda,\lambda+6}(t)=0.
\end{gather*}
We get the same results for $\lambda\in\{ a_1,
a_2,-\frac{5}{2},a_1-2,a_2-2,-\frac{9}{2},0,-2,-4,-3,-5,-7\}$ by
consi\-de\-ring the second-order integrability conditions.
\end{proof}

\begin{proposition}\label{pr6}For $k=9$, we have the following
third-order integrability  conditions of the infinitesimal
deformation~\eqref{InDefGen2}, for all $\lambda$:
\begin{alignat}{3}
& t_{\lambda,\lambda+3}\omega_{\lambda+3,\lambda+9}(t)=
\omega_{\lambda,\lambda+6}(t) t_{\lambda+6,\lambda+9}=0,& \nonumber\\
& t_{\lambda,\lambda+4} \omega_{\lambda+4,\lambda+9}(t) =
 \omega_{\lambda,\lambda+5}(t) t_{\lambda+5,\lambda+9} = 0,&\nonumber\\
& t_{\lambda-2,\lambda} \omega_{\lambda,\lambda+7}(t)=
 \omega_{\lambda,\lambda+7}(t) t_{\lambda+7,\lambda+9} = 0.&\label{thirdk9}
\end{alignat}
\end{proposition}

\begin{proof}
For $k=9$ and
$\lambda\notin\{0,-2,-4,-6,-8,a_i,-\frac{5}{2},a_i-3,-\frac{11}{2}\}$
we have
\begin{gather*}E_{\lambda,\lambda+9}= t_{\lambda,\lambda+3}
\omega_{\lambda+3,\lambda+9}(t)[\![J_{7}^{-1,\lambda+3},
C_{\lambda,\lambda+3}]\!]+t_{\lambda,\lambda+4}
\omega_{\lambda+4,\lambda+9}(t)[\![J_{6}^{-1,\lambda+4}, C_{\lambda,\lambda+4}]\!]\\
\phantom{E_{\lambda,\lambda+9}=}{} +t_{\lambda+5,\lambda+9}
\omega_{\lambda,\lambda+5}(t)[\![ C_{\lambda+5,\lambda+9},
J_{6}^{-1,\lambda}]\!]+t_{\lambda+6,\lambda+9}
\omega_{\lambda,\lambda+6}(t)[\![C_{\lambda+6,\lambda+9},
J_{7}^{-1,\lambda}]\!].
\end{gather*}
The equation $E_{\lambda,\lambda+9}=\partial J_{10}^{-1,\lambda}$
gives
\begin{gather*}
t_{\lambda,\lambda+3}
\omega_{\lambda+3,\lambda+9}(t)=t_{\lambda,\lambda+4}
\omega_{\lambda+4,\lambda+9}(t)=t_{\lambda+5,\lambda+9}
\omega_{\lambda,\lambda+5}(t)=t_{\lambda+6,\lambda+9}
\omega_{\lambda,\lambda+6}(t)=0.
\end{gather*}
By considering  the second-order integrability conditions, we get
the same results for each
$\lambda\in\{0,-2,-4,-6,-8,a_i,-\frac{5}{2},a_i-3,-\frac{11}{2}\}$.
\end{proof}

\begin{proposition}\label{pr7}For $k=10$, we have the following
third-order integrability  conditions of the infinitesimal
deformation~\eqref{InDefGen2}
\begin{alignat}{3}
& t_{\lambda-2,\lambda} \omega_{\lambda,\lambda+8}(t) =
 \omega_{\lambda,\lambda+8}(t) t_{\lambda+8,\lambda+10} = 0 \qquad&& \text{for}\quad
\lambda=0, -7,&\nonumber\\
& t_{\lambda-3,\lambda} \omega_{\lambda,\lambda+7}(t) =
\omega_{\lambda,\lambda+7}(t) t_{\lambda+7,\lambda+10} = 0
\qquad&& \text{for all}\quad \lambda,& \nonumber\\
& t_{\lambda,\lambda+4} \omega_{\lambda+4,\lambda+10}(t) =
 \omega_{\lambda,\lambda+6}(t) t_{\lambda+6,\lambda+10} = 0 \qquad&& \text{for all}\quad \lambda,&\nonumber\\
&t_{\lambda,\lambda+5} \omega_{\lambda+5,\lambda+10}(t) =
 \omega_{\lambda-5,\lambda}(t)\,t_{\lambda,\lambda+5}= 0 \qquad&& \text{for}\quad
\lambda=0,-4.&\label{thirdk10}
\end{alignat}
\end{proposition}
\begin{proof}
For $k=10$ and
$\lambda\notin\{-9,-7,-\frac{13}{2},-6,-5,-4,-3,-2,-\frac{5}{2},0,a_1,a_1-4,a_2,a_2-4\}$
we have
\[
E_{\lambda,\lambda+10}=t_{\lambda,\lambda+4}
\omega_{\lambda+4,\lambda+10}(t)[\![J_{7}^{-1,\lambda+4},
C_{\lambda,\lambda+4}]\!]+t_{\lambda+6,\lambda+10}
\omega_{\lambda,\lambda+6}(t)
[\![C_{\lambda+6,\lambda+10},J_{7}^{-1,\lambda}]\!].
\]
The equation $E_{\lambda,\lambda+10}=\partial J_{11}^{-1,\lambda}$
gives the conditions
\[
t_{\lambda,\lambda+4}
\omega_{\lambda+4,\lambda+10}(t)=t_{\lambda+6,\lambda+10}
\omega_{\lambda,\lambda+6}(t)=0.
\]
We check that, for
$\lambda\in\{-9,-7,-\frac{13}{2},-6,-5,-4,-3,-2,-\frac{5}{2},0,a_1,a_1-4,a_2,a_2-4\},$
these latter conditions must be also satisf\/ied. The others
conditions follow from the singular values of~$\lambda$.
\end{proof}

\begin{proposition}\label{pr8}For $k=11$, we have the following
second-order integrability  conditions of the infinitesimal
deformation~\eqref{InDefGen2}
\begin{alignat}{3}
& t_{\lambda,\lambda+6}\,\omega_{\lambda+6,\lambda+11}(t) =
\omega_{\lambda-5,\lambda}(t) t_{\lambda,\lambda+6} = 0 \qquad&& \text{for}\quad \lambda=-\tfrac{5\pm \sqrt{19}}{2},&\nonumber\\
& t_{\lambda,\lambda+5}\,\omega_{\lambda+5,\lambda+11}(t) =
 \omega_{\lambda-6,\lambda}(t) t_{\lambda,\lambda+5} = 0 \qquad && \text{for}\quad \lambda=0, -4,&\nonumber\\
& t_{\lambda-4,\lambda} \omega_{\lambda,\lambda+7}(t) =
 \omega_{\lambda,\lambda+7}(t) t_{\lambda+7,\lambda+11} = 0
\qquad&& \text{for all}\quad \lambda,& \nonumber\\
& t_{\lambda-3,\lambda} \omega_{\lambda,\lambda+8}(t) =
 \omega_{\lambda,\lambda+8}(t) t_{\lambda+8,\lambda+11} = 0 \qquad && \text{for}\quad
\lambda=0, -7.& \label{thirdk11}
\end{alignat}
\end{proposition}

\begin{proposition}\label{pr9}For $k=12$, we have the following
third-order integrability  conditions of the infinitesimal
deformation~\eqref{InDefGen2}
\begin{alignat}{3}
&t_{\lambda,\lambda+5}\,\omega_{\lambda+5,\lambda+12}(t)=
\omega_{\lambda-7,\lambda}(t)\,t_{\lambda,\lambda+5}=0,
\qquad&&\text{for}\quad \lambda=0,\,-4,&\nonumber\\
& t_{\lambda,\lambda+6} \omega_{\lambda+6,\lambda+12}(t) =
\omega_{\lambda-6,\lambda}(t) t_{\lambda,\lambda+6} = 0
\qquad&& \text{for}\quad \lambda=-\tfrac{5\pm \sqrt{19}}{2},&\nonumber\\
& t_{\lambda-4,\lambda} \omega_{\lambda,\lambda+8}(t) =
 \omega_{\lambda,\lambda+8}(t) t_{\lambda+8,\lambda+12} = 0 \qquad&& \text{for}\quad
\lambda=0,-7.&\label{k12}
\end{alignat}
\end{proposition}

\begin{proof}
 For $k=11$ and $k=12$, the 2-cochains
$E_{\lambda,\lambda+k}$ are def\/ined only for some particular
values of $\lambda$. We compute these 2-cochains and then we check
the corresponding integrability conditions.
\end{proof}
\subsection{Fourth-order integrability conditions}

\begin{proposition}\label{pr10} For generic $\lambda$, the fourth-order integrability  conditions of the infinitesimal
deformation~\eqref{InDefGen2} are the following:
\begin{equation}\label{fourth}
\omega_{\lambda,\lambda+i}(t) \omega_{\lambda+i,\lambda+k}(t)=0,\qquad\text{where}\quad 5\leq
i\leq7\qquad\text{and}\qquad 5+i\leq k\leq7+i.
 \end{equation}
\end{proposition}
\begin{proof}
These conditions come from the fact that the fourth term
$\mathcal{L}^{(4)}$ must satisfy:
\[
\partial {\cal L}^{(4)}=-\frac{1}{2}\,[\![{\cal
L}^{(2)},{\cal L}^{(2)}]\!].
\]
Indeed, we can always reduce $\mathcal{L}^{(3)}$ to zero by
equivalence.
\end{proof}
The following theorem is our main result.
\begin{theorem}
\label{Main2}
 The second-order integrability conditions \eqref{k567} and \eqref{k89}
 together with the third and the fourth-order
conditions \eqref{thirdk7}--\eqref{fourth} are necessary and
sufficient for the complete integrability of the infinitesimal
deformation~\eqref{InDefGen2}. Moreover, any formal ${\rm
\frak{sl}}(2)$-trivial deformation of the Lie derivative~$L_X$ on
the space of symbols ${\cal S}_\delta^n$ is equivalent to a
polynomial one of degree equal or less than~$2$.
\end{theorem}

\begin{proof} Clearly, all these conditions are necessary.
So, let us prove that they are also suf\/f\/icient. As in the
proof of Theorem \ref{Main1}, the solution $\mathcal{L}^{(3)}$ of
the Maurer--Cartan equation (\ref{MC2}) is def\/ined up to a
1-coboundary, thus, we can always reduce $\mathcal{L}^{(3)}$ to
zero by equivalence. Moreover, by recurrence, the highest-order
terms $\mathcal{L}^{(m)}$ satisfy the equation
$\partial\mathcal{L}^{(m)}=0$ and can also be reduced to the
identically zero map. This completes the proof of
Theorem~\ref{Main2}.
\end{proof}

\begin{remark}\label{remark} The majority of integrability conditions
concern some parameters $t_{\lambda,\lambda+k}$ with singular
values of $\lambda$. All these singular values of $\lambda$ are
negatives. So, let us consider the space ${\cal S}^n_{\delta}$
with generic $\delta$, for example, $\delta-n>0$. In this case, the second-order integrability conditions are reduced
to the following equations:
\[
b t_{\lambda,\lambda+3}t_{\lambda+3,\lambda+7}-a t_{\lambda,\lambda+4}t_{\lambda+4,\lambda+7}=
t_{\lambda,\lambda+4}t_{\lambda+4,\lambda+8}=0.
\]
\end{remark}

\section{Examples}\label{sec7}

\begin{example} \label{Example1} Let us consider  the space of symbols
${\cal S}^4_{\lambda+4}$.

\begin{proposition} Any formal ${\rm \frak{sl}}(2)$-trivial
deformation of the ${\rm Vect_{Pol}}(\mathbb{R})$-action on the
spa\-ce~${\cal S}^4_{\lambda+4}$ is equivalent to his
infinitesimal part, without any conditions on the parameters
(independent parameters). That is, the miniversal deformation
here has base $\mathbb{C}[[t]]$ where $t$ designates the family
of all parameters.
\end{proposition}

\begin{proof} The inf\/initesimal ${\rm \frak{sl}}(2)$-trivial deformation, in this case, is given by
\[
{\cal L}=L+{\cal L}^{(1)},
\]
where $L_X$ is the Lie derivative of ${\cal S}^4_{\lambda+4}$
along the vector f\/ield $X\frac{d}{dx}$ def\/ined
by~(\ref{Lieder2}), and
\begin{gather*}
{\cal L}^{(1)}=t_{\lambda,\lambda+2} C_{\lambda,\lambda+2}+
t_{\lambda,\lambda+3} C_{\lambda,\lambda+3}+t_{\lambda,\lambda+4} C_{\lambda,\lambda+4}\nonumber\\
\phantom{{\cal L}^{(1)}=}{} +t_{\lambda+1,\lambda+3}
C_{\lambda+1,\lambda+3}(X)+
t_{\lambda+1,\lambda+4} C_{\lambda+1,\lambda+4}+t_{\lambda+2,\lambda+4} C_{\lambda+2,\lambda+4}.
\end{gather*}
There are no conditions to integrate this inf\/initesimal
deformation to a formal one. The solution~$\mathcal{L}^{(2)}$
of~(\ref{MCG}) is def\/ined up to a 1-coboundary and dif\/ferent
choices of solutions of the Maurer--Cartan equation correspond to
equivalent deformations. Thus, we can always
reduce~$\mathcal{L}^{(2)}$ to zero by equivalence. Then, by
recurrence, the highest-order terms $\mathcal{L}^{(m)}$ satisfy
the equation $\partial\mathcal{L}^{(m)}=0$ and $\mathcal{L}^{(m)}$
can also be reduced to the identically zero map.
\end{proof}
\end{example}

\begin{remark} We have the same results for ${\cal S}^k_{\lambda+k}$ if
$k\leq4.$ Indeed, for $k\leq4$ there are no integrability
conditions.
\end{remark}

\begin{example} \label{Example2} Let us consider  the ${\rm
Vect_{Pol}}(\mathbb{R})$-module ${\cal S}^6_{7}$.

\begin{proposition} Any formal ${\rm \frak{sl}}(2)$-trivial
deformation of the ${\rm Vect_{Pol}}(\mathbb{R})$-action on the
space ${\cal S}^6_{7}$ is equivalent to
\begin{gather}\label{l67}
{\cal L}_X=L+{\cal L}^{(1)}+{\cal L}^{(2)},
\end{gather}
where \begin{gather*} {\cal L}^{(1)}= t_{1,3} C_{1,3}+
t_{1,4} C_{1,4}+t_{1,5} C_{1,5}+t_{2,4} C_{2,4}+ t_{2,5} C_{2,5}+t_{2,6} C_{2,6}\\
\phantom{{\cal L}^{(1)}=}{} +t_{3,5} C_{3,5}+ t_{3,6}
C_{3,6}+t_{3,7} C_{3,7}+ t_{4,6} C_{4,6}+t_{4,7} C_{4,7}
\end{gather*}
 and
\begin{gather*}{\cal
L}^{(2)}
=\tfrac{1}{2}\omega_{1,6}(t)J_6^{-1,1}+\tfrac{1}{2}\omega_{2,7}(t)J_6^{-1,2}+
\tfrac{1}{2}\omega_{1,7}(t)J_7^{-1,1}.
\end{gather*}

The formal deformation \eqref{l67} is defined without any
condition on the parameters (independent parameters). That is, the
miniversal deformation  here has base $\mathbb{C}[[t]]$ where
$t$ designates the family of all parameters.
\end{proposition}
\end{example}

For $k\leq6$, generically there are no integrability conditions
which is the case of the previous example (see
Remark~\ref{remark}). Now, we study a generic example with $k=7$.

\begin{example} \label{Example3} Let us consider  the ${\rm Vect_{Pol}}(\mathbb{R})$-module  ${\cal
S}^7_{\lambda+7}$ for generic $\lambda$.

\begin{proposition}The ${\rm Vect_{Pol}}(\mathbb{R})$-module
 ${\cal S}^7_{\lambda+7}$
admits six ${\rm \frak{sl}}(2)$-trivial deformations with
$11$~independent parameters. It admits a miniversal ${\rm
\frak{sl}}(2)$-trivial deformation with $15$  parameters. These
deformations are polynomial of degree~$2$.
\end{proposition}

\begin{proof} Any formal ${\rm \frak{sl}}(2)$-trivial deformation of the
${\rm Vect_{Pol}}(\mathbb{R})$-action on the space ${\cal
S}^7_{\lambda+7}$ is equiva\-lent to
\begin{gather*}
{\cal L}=L+{\cal L}^{(1)}+{\cal L}^{(2)},
\end{gather*}
where \begin{gather*}  {\cal
L}^{(1)}=t_{\lambda,\lambda+2}C_{\lambda,\lambda+2}+
t_{\lambda,\lambda+3}C_{\lambda,\lambda+3}+t_{\lambda,\lambda+4}C_{\lambda,\lambda+4}
+t_{\lambda+1,\lambda+3}C_{\lambda+1,\lambda+3}\nonumber\\
 \phantom{{\cal L}^{(1)}=}{} +t_{\lambda+1,\lambda+4} C_{\lambda+1,\lambda+4}+t_{\lambda+1,\lambda+5} C_{\lambda+1,\lambda+5}
+t_{\lambda+2,\lambda+4} C_{\lambda+2,\lambda+4}+
t_{\lambda+2,\lambda+5} C_{\lambda+2,\lambda+5}\nonumber\\
\phantom{{\cal L}^{(1)}=}{}+ t_{\lambda+2,\lambda+6}
C_{\lambda+2,\lambda+6}+t_{\lambda+3,\lambda+5}
C_{\lambda+3,\lambda+5} +t_{\lambda+3,\lambda+6}
C_{\lambda+3,\lambda+6}+
t_{\lambda+3,\lambda+7} C_{\lambda+3,\lambda+7}\nonumber\\
\phantom{{\cal L}^{(1)}=}{}+t_{\lambda+4,\lambda+6}
C_{\lambda+4,\lambda+6} +t_{\lambda+4,\lambda+7}
C_{\lambda+4,\lambda+7}
+t_{\lambda+5,\lambda+7} C_{\lambda+5,\lambda+7} 
\end{gather*}
 and
\begin{gather*}
{\cal
L}^{(2)}=\tfrac{1}{2}\omega_{\lambda,\lambda+5}(t)J_6^{-1,\lambda}+
\tfrac{1}{2}\omega_{\lambda+1,\lambda+6}(t)J_6^{-1,\lambda+1}+
\tfrac{1}{2}\omega_{\lambda+2,\lambda+7}(t)J_6^{-1,\lambda+2}\\
\phantom{{\cal L}^{(2)}=}{} +
\tfrac{1}{2}\omega_{\lambda,\lambda+6}(t)J_7^{-1,\lambda}+
\tfrac{1}{2}\omega_{\lambda+1,\lambda+7}(t)J_7^{-1,\lambda+1}.
\end{gather*}
There are only three integrability conditions:
\begin{gather*}
b\,t_{\lambda,\lambda+3}t_{\lambda+3,\lambda+7}-a\,t_{\lambda,\lambda+4}t_{\lambda+4,\lambda+7}=t_{\lambda,\lambda+2}\,\omega_{\lambda+2,\lambda+7}(t)=\omega_{\lambda,\lambda+5}(t)\,t_{\lambda+5,\lambda+7}=0.
\end{gather*}
The formal ${\rm \frak{sl}}(2)$-trivial deformations with the
greatest number of independent parameters are those corresponding
to
$b\,t_{\lambda,\lambda+3}t_{\lambda+3,\lambda+7}-a\,t_{\lambda,\lambda+4}t_{\lambda+4,\lambda+7}=t_{\lambda,
\lambda+2}=t_{\lambda+5, \lambda+7}=0$.  So, we must kill at least
four parameters and there are six choices. Thus, there are  four
deformations with 11~independent parameters. Of course, there are
many formal deformations with less then 11~independent parameters.
The deformation ${\cal L}=L+{\cal L}^{(1)}+{\cal L}^{(2)}, $ is
the miniversal $\mathfrak{sl}(2)$-trivial deformation of ${\cal
S}^{7}_{\lambda+7}$ with base $\mathcal{A}=
\mathbb{C}[t]/\mathcal{R}$, where $t$ is the family of all
parameters given in the expression of ${\cal L}^{(1)}$ and
$\mathcal{R}$ is the ideal generated by the polynomials
$b\,t_{\lambda,\lambda+3}t_{\lambda+3,\lambda+7}-a\,t_{\lambda,\lambda+4}t_{\lambda+4,\lambda+7}$ and $
\omega_{\lambda,\lambda+5}(t) t_{\lambda+5,\lambda+7}$.
\end{proof}
\end{example}

\subsection*{Acknowledgements} We would like to thank Sofiane Bouarroudj for his suggestions and for pointing out a mistake in a previous version of this paper.

\pdfbookmark[1]{References}{ref}
\LastPageEnding


\begin{thebibliography}{99}

\footnotesize\itemsep=0pt

\bibitem{abbo} Agrebaoui B., Ben Fraj N., Ben Ammar M.,   Ovsienko  V.,
Deformation of modules of dif\/ferential forms, {\it J.~Nonlinear
Math. Phys.} {\bf 10} (2003), 148--156,
\href{http://arxiv.org/abs/math.QA/0310494}{math.QA/0310494}.

\bibitem{aalo2}
 Agrebaoui B., Ammar F., Lecomte P.,  Ovsienko V.,
Multi-parameter deformations of the module of symbols of
dif\/ferential operators, {\it Int. Math. Res. Not.} {\bf 2002}
(2002), no.~16, 847--869,
\href{http://arxiv.org/abs/math.QA/0011048}{math.QA/0011048}.


\bibitem{b}
 Bouarroudj S., On $\frak {sl}(2)$-relative cohomology of the Lie
algebra of vector f\/ields and dif\/ferential operators, {\it
J.~Nonlinear Math. Phys.} {\bf 14} (2007), 112--127,
\href{http://arxiv.org/abs/math.DG/0502372}{math.DG/0502372}.

\bibitem{bo}
 Bouarroudj S.,  Ovsienko V.,  Three cocycles on Dif\/f($S^1$)
generalizing the Schwarzian derivative, {\it Int. Math. Res. Not.}
{\bf 1998} (1998), no.~1, 25--39,
\href{http://arxiv.org/abs/dg-ga/9710018}{dg-ga/9710018}.


\bibitem{f1}
 Fialowski A.,  Deformations of Lie algebras, {\it Mat. Sb.}
{\bf 55} (1986), 467--473.

\bibitem{f2}
 Fialowski A.,  An example of formal deformations of Lie
algebras, in  Deformation Theory of Algebras and Structures and
Applications, {\it NATO Adv. Sci. Inst. Ser. C Math. Phys. Sci.},
Vol.~247, Kluwer Acad. Publ., Dordrecht, 1988, 375--401.

\bibitem{ff2}
 Fialowski A., Fuchs D.B., Construction of miniversal
deformations of Lie algebras, {\it J. Funct. Anal.} {\bf 161}
(1999), 76--110,
\href{http://arxiv.org/abs/math.RT/0006117}{math.RT/0006117}.

\bibitem{Fuc2}
Fuchs D.B., Cohomology of inf\/inite-dimensional Lie algebras,
Consultants Bureau, New York, 1987.

\bibitem{hg}
 Gargoubi H., Sur la g\'eom\'etrie de l'espace des
op\'erateurs dif\/f\'erentiels lin\'eaires sur $\mathbb{R}$, {\it
Bull. Soc. Roy. Sci. Li\`ege} {\bf 69} (2000), 21--47.

\bibitem{gmo}
Gargoubi H.,  Mellouli N., Ovsienko V., Dif\/ferential operators
on supercircle: conformally equivariant quantization and symbol
calculus, {\it Lett. Math. Phys.} {\bf 79} (2007), 51--65,
\href{http://arxiv.org/abs/math-ph/0610059}{math-ph/0610059}.

\bibitem{pg}
 Gordan P., Invariantentheorie, Teubner, Leipzig, 1887.

\bibitem{nr2}
 Nijenuis A., Richardson  R.W. Jr., Deformations of
homomorphisms of Lie groups and Lie algebras, {\it Bull. Amer.
Math. Soc.} {\bf 73} (1967), 175--179.

\bibitem{ro12}
Ovsienko V., Roger C., Deforming the Lie algebra of vector
f\/ields on~$S^1$ inside the Lie algebra of pseudodif\/ferential
operators on $S^1$, in Dif\/ferential Topology,
Inf\/inite-Dimensional Lie Algebras, and Applications, {\it Amer.
Math. Soc. Transl. Ser. 2}, Vol.~194, Amer. Math. Soc.,
Providence, RI, 1999, 211--226,
\href{http://arxiv.org/abs/math.QA/9812074}{math.QA/9812074}.

\bibitem{ro22}
Ovsienko V., Roger C., Deforming the Lie algebra of vector
f\/ields on~$S^1$ inside the Poisson algebra on $\dot T^* S^1$,
{\it Comm. Math. Phys.} {\bf 198} (1998), 97--110,
\href{http://arxiv.org/abs/q-alg/9707007}{q-alg/9707007}.


\bibitem{r}
Richardson R.W., Deformations of subalgebras of Lie algebras, {\it
J. Differential Geom.} {\bf 3} (1969), 289--308.

\end{thebibliography}
\end{document}